%% file: lotoreichik_ourmieres.tex
\documentclass[a4paper,11pt]{amsart}
\usepackage{amssymb,amsmath,amsthm}
\usepackage{fullpage}
\usepackage[colorlinks=true,linkcolor=blue, citecolor=red]{hyperref}
\usepackage{eucal}
\usepackage{color}
\usepackage{xcolor}
\usepackage{stackrel, enumerate, mathrsfs, lmodern}
\usepackage[normalem]{ulem}
\usepackage{graphicx, bbm}
\usepackage[left=2.7cm, right=2.7cm, marginparwidth=2cm, footskip = 0.8cm, textheight=700pt ]{geometry}
\definecolor{darkblue}{rgb}{0.2,0.2,0.6}
\definecolor{superdarkblue}{rgb}{0.2,0.2,0.3}
\hypersetup{colorlinks=true, linkcolor=darkblue,citecolor=darkblue, urlcolor = superdarkblue}
\definecolor{citegreen}{rgb}{0.2,0.2,0.6}
\definecolor{green2}{rgb}{0.,0.6,0.2}

\input commands.tex

\newcommand{\Hm}[1]{\leavevmode{\marginpar{\scriptsize%
			$\hbox to 0mm{\hspace*{-0.5mm}$\leftarrow$\hss}%
			\vcenter{\vrule depth 0.1mm height 0.1mm width \the\marginparwidth}%
			\hbox to
			0mm{\hss$\rightarrow$\hspace*{-0.5mm}}$\\
			\relax\raggedright #1}}}

\newcommand\Th{\Theta}

\newtheorem{corol}{Corollary}[section]

\title[]{Spectral asymptotics of the Dirac operator in a thin shell}
\author{Vladimir Lotoreichik}

\address[V. Lotoreichik]{Department of Theoretical Physics, Nuclear Physics Institute, Czech Academy of Sciences,
25068 \v{R}e\v{z}, Czech Republic}
\email{lotoreichik@ujf.cas.cz}
\urladdr{http://gemma.ujf.cas.cz/~lotoreichik/}
\author{Thomas Ourmi\`{e}res-Bonafos}
\address[T. Ourmi\`eres-Bonafos]{Aix-Marseille Universit\'e, CNRS, Centrale Marseille, I2M, Marseille, France.}
\email{thomas.ourmieres-bonafos@univ-amu.fr}
\urladdr{http://www.i2m.univ-amu.fr/perso/thomas.ourmieres-bonafos/}

\begin{document}
\maketitle
\begin{abstract} We investigate the spectrum of the Dirac operator with infinite mass boundary conditions posed in a tubular neighborhood of a smooth compact hypersurface in $\RR^n$ without boundary.  We prove that when the tubular neighborhood shrinks to the hypersurface, the asymptotic behavior of the eigenvalues is driven by a Schr\"odinger operator involving electric and Yang-Mills potentials, both of geometric nature.
\end{abstract}
\tableofcontents
\section{Introduction and main result}
\subsection{Motivations and state of the art}
Because of possible applications in mesoscopic physics, quantum systems in thin structures have been extensively studied in the last thirty years and a canonical example is the Dirichlet Laplacian posed in a tubular neighborhood of a curve in dimension two or three. Such a system is usually called a quantum waveguide because this operator is the Hamiltonian of a spinless quantum particle trapped in the waveguide formed by this tubular neighborhood. As this particle obeys the Schr\"odinger equation, it is well known that its dynamics is given by the structure of the spectrum of this Hamiltonian and in the works \cite{ES89,GJ92}, it is shown that the curvature of the underlying curve plays a crucial role in the apparition of bound states. In \cite{DE95}, this study is pushed further and it is shown that when the waveguide degenerates onto the underlying curve an effective Schr\"odinger operator of geometric nature drives the dynamics.  Since then, similar questions have been addressed in quantum layers, that is in the tubular neighborhood of a hypersurface in $\RR^n$. Again, as shown in \cite{DEK01,EK01,CEK04}, the geometry of the hypersurface can create bound states and its curvatures play an important role. Further results and references on quantum waveguides and layers
can be found in the monograph~\cite{EK15}.

This is a similar geometrical setting we investigate in this paper but instead of the Dirichlet Laplacian, we consider the Dirac operator with the so-called infinite mass boundary conditions posed in the tubular neighborhood of a smooth compact hypersurface in $\RR^n$ without boundary. Our aim is to highlight the influence of the geometry of the underlying hypersurface on the spectrum of this operator when the tubular neighborhood degenerates onto the hypersurface.

From a physical point of view, the Dirac operator models a relativistic quantum particle of spin-$\frac12$ and this operator plays an important role in the study of two-dimensional structures having honeycomb symmetries such as graphene (see, for instance \cite{FW12}) and is used as a phenomenological model for the confinement of quarks in hadrons in dimension three (see, \cite{Bogo,CJKTW}). About the boundary conditions we consider, it has been remarked in \cite{BM87} that they arise when considering the Dirac operator posed in the full Euclidean space $\RR^n$ and letting a mass supported outside the considered domain tend to infinity. This has been rigorously justified in the series of papers \cite{SV19,BCLTS19,ALTMR19}. In three dimensions these boundary conditions appear  in the MIT bag model describing the confinement of quarks~\cite{J75}, while in two dimensions the same type of boundary conditions arises in the description of graphene quantum dots (see~\cite{BFSV17} and the references therein).  Regarding asymptotic regimes exhibiting geometrical properties for the spectrum of Dirac operators we mention the seminal paper \cite{ALTR17} in which the authors consider, in dimension three, a regime where the boundary of the domain becomes attractive. This work is generalized in any dimension in \cite{MOP20} and in the abstract framework of spin-manifolds in \cite{Fla22}. We also mention \cite{HOP18} where a highly attractive zero-range potential is considered in dimension three as well as \cite{BLTRS23} where strong magnetic fields in dimension two are investigated.

To our knowledge, the only cases where the spectrum of the Dirac operator with infinite mass boundary conditions is studied in the regime of a tubular neighborhood shrinking to the underlying hypersurface is for two-dimensional waveguides in \cite{BBKOB22,BKOB23}. It turns out that a non-commutative gauge transform allows to simplify the geometrical properties of the effective operator which is no longer the case when the problem is set in the tubular neighborhood of a closed loop as shown in \cite{K23}. This is the precise question we generalize in this paper to any space dimension.

\subsection{The Dirac operator in a thin shell}\label{subsec:dots}
We start by defining the geometry we are interested in. Let $\Sg\subset\dR^n$, $n\ge 2$, be a $C^\infty$-smooth compact hypersurface without boundary. Let $\nu$ be the normal vector field to $\Sigma$ pointing outward the bounded domain delimited by $\Sigma$. The shape operator $S : T\Sigma \to T\Sigma$ is given by $SX = \nabla_X \nu$ and its eigenvalues $\kappa_1,\dots,\kappa_{n-1}$ are the principal curvatures of $\Sigma$. For $p \in \{1,\dots,n-1\}$, $H_p$ denotes the $p$-th mean curvature of $\Sigma$ with respect to $\nu$ and verifies
\[
	H_p = \sum_{1\leq j_1<\cdots<j_p \leq n-1}\kappa_{j_1}\cdot \cdots \cdot\kappa_{j_p}.
\]
By convention, we set $H_p = 0$ for $p \geq n$. 

For $\eps > 0$ consider the map
\begin{equation}
	\phi_\eps\colon \Sg\times(-1,1)\arr\dR^n,\qquad 
	\phi_\eps(s,t) := s +\eps t\nu(s),
	\label{eqn:tubcoord}
\end{equation}
where $\nu(s)$ is the unit normal vector to $\Sg$ at $s\in\Sg$ pointing outwards of the bounded domain surrounded by $\Sg$. According to~\cite[Theorem 10.19]{Lee} there exists $\eps_0 > 0$ such that for all $\eps\in(0,\eps_0)$ the map $\phi_\eps$ is a $C^\infty$-diffeomorphism from $\Sg\times(-1,1)$ onto the shell $\Omg_\eps := \phi_\eps(\Sg\times(-1,1))$.

We are interested in the Dirac operator posed in the shell $\Omega_\eps$ with infinite mass boundary conditions on $\partial\Omega_\eps$. To introduce it, set $N = 2^{\lfloor\frac{n+1}{2}\rfloor}$ and let $(\aa_j)_{j\in\{1,\dots,n+1\}}$ be hermitian matrices of size $N\times N$ verifying for all $j,k\in\{1,\dots,n+1\}$:
\begin{equation}
	\alpha_j\alpha_k + \alpha_k\alpha_j = 2 \delta_{j,k},
	\label{eqn:anticomalpha}
\end{equation}
where $\delta_{j,k}$ is the standard Kronecker symbol.
Such matrices are constructed explicitly further in \S \ref{subsec:eqnDiracmat}.

For any $x\in\dR^n$ we define the matrix $\G(x) := \sum_{j=1}^n x_j\aa_j$ and for $\eps \in(0,\eps_0)$, we introduce the Dirac operator acting in $L^2(\Omg_\eps;\dC^N)$ by
\begin{equation}\label{eqn:defdirac}
\begin{aligned}
 \scD_\eps u& := -i\sum_{j=1}^n\aa_j\p_j u + m\aa_{n+1}u,\\
 \dom\scD_\eps &:= \left\{ u\in H^1(\Omg_\eps;\dC^N)\colon -i\aa_{n+1}\G(\nu_\eps)u = u\,\,\,\text{on}\,\,\p\Omg_\eps\right\},
\end{aligned}
\end{equation}
where $\nu_\eps$ is the outer unit normal vector to $\Omg_\eps$ and the mass $m \geq 0$ is a fixed parameter.

Thanks to \cite[Ex. 4.20 \& Thm. 4.11]{BB16}, $\scD_\eps$ is known to be self-adjoint and it has compact resolvent. Hence, its spectrum is constituted of discrete eigenvalues and the purpose of this paper is to understand the behavior of these eigenvalues in the thin shell regime $\eps \to 0$. It turns out that in this regime an effective operator on $\Sigma$ plays an important role. Namely, consider the matrix-valued $1$-form on $\Sigma$ given by the local expression
\begin{equation}\label{eq:Omega}
	\Omega = \sum_{j=1}^{n-1}\Omega^j \dd s_j \in T^*\Sigma\otimes \CC^{N\times N},\qquad \Omega^j = -i \Gamma(\partial_j\nu) \Gamma(\nu)\in\dC^{N\times N}
\end{equation}
and the Schr\"odinger operator with an external Yang-Mills potential $\Upsilon$ defined by Kato's first representation theorem (see \cite[Ch. VI, Thm. 2.1]{Kat}) as the unique self-adjoint operator associated with the closed sesquilinear form given by
\begin{equation}
\label{eq:theta}
\begin{aligned}
	\mathfrak{u}[f,g] &:= \int_\Sigma \sum_{j,k=1}^{n-1} g^{j,k}\left\langle\partial_j f - i \left(\frac12 - \frac1\pi\right)\Omega^jf,\partial_kg - i \left(\frac12 - \frac1\pi\right)\Omega^kg\right\rangle_{\CC^N}\dd s \\&\qquad\qquad\qquad\qquad\qquad+ \int_\Sg\left(\left(\frac12 + \frac2{\pi^2}\right)H_2 - \frac{H_1^2}{\pi^2}\right) \langle f,g\rangle_{\CC^N} \dd s,\\ \dom{\mathfrak{u}} &:= H^1(\Sigma, \CC^N).
\end{aligned}
\end{equation}
Here $g = (g_{j,k})_{j,k \in \{1,\dots,n-1\}}$ is the Riemannian metric on $\Sigma$ induced by the embedding into $\RR^n$ and $(g^{j,k})_{j,k \in \{1,\dots,n-1\}} := g^{-1}$.

To state our main result, we recall the notion of min-max levels: for a closed quadratic form $Q$, semi-bounded below with dense domain $\dom{Q}$ in a Hilbert space $\cH$ we define, for $p\in \N$, the $p$-th min-max level of $Q$ as
\begin{equation}
	\mu_p(Q) := \inf_{\substack{F \subset \dom{Q}\\\dim F = p}} \sup_{u \in F\setminus\{0\}} \frac{Q[u]}{\|u\|_{\cH}^{2}}.
	\label{eqn:minmax}
\end{equation}
By Kato's first representation theorem, the sesquilinear form associated with $Q$ gives rise to a unique self-adjoint operator $A$ acting in $\cH$ (see \cite[Ch. VI, Thm. 2.1]{Kat}). If $A$ has compact resolvent, for $p\in \N$, $\mu_p(Q)$ is the $p$-th eigenvalue of $A$ counted with multiplicities. In this case, we set
\begin{equation}
	\mu_p(A) := \mu_p(Q).
	\label{eqn:notminmax}
\end{equation}

Our main result is about the asymptotic expansion of the eigenvalues of the operator $\scD_\eps^2$. It reads as follows.
\begin{thm} For all $p\in \N$, there holds:
\[
	\mu_p(\scD_\eps^2) = \frac{\pi^2}{16\eps^2} + \frac{m}\eps -\frac{4}{\pi^2} m^2 +m^2 + \mu_p(\Upsilon) + \mathcal{O}(\varepsilon),\qquad \eps \to 0.
\]
\label{thm:main}
\end{thm}
When $\eps \to 0$ the domain $\Omega_\eps$ collapses to the hypersurface $\Sigma$ and one of the main features of Theorem \ref{thm:main} is that up to a renormalization factor the spectrum of $\scD_\eps^2$ behaves like the one of the operator $\Upsilon$ which is, by definition, an object involving geometric features of $\Sigma$. In particular, the splitting of the eigenvalues appears thanks to the effective operator $\Upsilon$ at order $0$ whereas the main term in the expansion is of order $\eps^{-2}$.

A natural question is to understand whether from the expansion of the eigenvalues of $\scD_\eps^2$ given in Theorem \ref{thm:main} one can obtain asymptotic expansions for the eigenvalues of $\scD_\eps$. As explained in \cite[Prop. A.2.]{MOP20}, if $n \notin 4\ZZ$, due to the existence of a so-called {\it charge conjugation} operator, the spectrum of $\scD_\eps$ is symmetric with respect to the origin and if $E \in Sp(\scD_\eps)$ then $\dim\big({\ker(\scD_\eps - E)}\big) = \dim\big({\ker(\scD_\eps + E)}\big)$. It yields the following corollary.
\begin{corol} Let $n \in \N$ be such that $n \notin 4 \ZZ$. Then, if for $p\in \N$, $\lambda_p(\scD_\eps)$ denotes the $p$-th non-negative eigenvalue of $\scD_\eps$ counted with multiplicities, there holds
\[
	\lambda_p(\scD_\eps) = \frac{\pi}{4\varepsilon} + \frac2{\pi}m + \Big(\frac{2}\pi \mu_{2p}(\Upsilon) + \frac2\pi m^2 - \frac{16}{\pi^3}m^2\Big)\varepsilon + \mathcal{O}(\varepsilon^2).
\]
\end{corol}
When $n \in 4 \ZZ$ there is no reason for a {\it charge conjugation} operator to exist (see \cite[Thm. 15.19]{DG13}). Hence the spectrum $Sp(\scD_\varepsilon)$ of $\scD_\varepsilon$ may not be symmetric with respect to the origin and our approach does not allow to recover the full structure of the spectrum of $\scD_\varepsilon$.\\

The proof of Theorem \ref{thm:main} is by obtaining an upper and a lower bound on the min-max levels of the quadratic form of the square of the operator $\scD_\eps$. The first step is to rewrite the problem in the tubular coordinates $(s,t)$ introduced by the map $\phi_\eps$ given in \eqref{eqn:tubcoord} and to perform a unitary transform to obtain a metric on $\Sigma\times(-1,1)$ of the form $\dd s^2 + \dd t^2$, where $\dd s^2$ is the Riemannian metric on $\Sigma$ induced by its embedding into $\RR^n$ (up to lower order terms in $\eps$). Once this is done, one realizes that the leading term in the asymptotic expansion of the min-max levels of $\scD_\eps^2$ is driven by a one-dimensional operator acting in the transverse variable $t$ only. Thus, a second step consists in a thorough description of the spectral properties of this transverse operator. Note that this step is already  technically demanding because it uses extensively the explicit structure of the Dirac matrices $(\alpha_j)_{j\in\{1,\dots,n+1\}}$ constructed further in \S \ref{subsec:eqnDiracmat}. The next step is to find an upper bound using a trial function chosen as a tensor product between the modes associated with the lowest eigenvalue of the transverse operator and functions in the $s$-variable only. This step requires for detailed computations in local coordinates in order to make appear the effective operator $\Upsilon$. The last step concerns the lower bound and our goal is attained when controlling terms which are in the orthogonal complement of the vector space spanned by the modes associated with the lowest eigenvalue of the transverse operator. The main difficulty in this step is that the transverse operator depends on the $s$-variable through its boundary conditions leading to commutator estimates, obtained working in local coordinates on the hypersurface $\Sigma$.
\subsection{Organization of the paper} In Section \ref{sec:prelimi} we gather preliminary material about the Dirac matrices $(\alpha_j)_{j\in \{1,\dots,n+1\}}$, standard results of differential geometry that are used further in the paper and we perform the spectral analysis of the transverse operator playing a central role in the first-order terms of the asymptotic expansion of Theorem \ref{thm:main}. Section \ref{sec:qf} deals with the quadratic form of the operator $\scD_\eps^2$ which we rewrite in tubular coordinates using the map $\phi_\eps$ defined in \eqref{eqn:tubcoord} and a unitary transform to obtain an adequate metric on $\Sigma\times (-1,1)$. We obtain an upper and a lower bound on the operator $\scD_\eps^2$ in the form sense.  Finally, these bounds are used to prove Theorem \ref{thm:main} in Section \ref{sec:proofmr}. The paper is complemented by Appendix~\ref{app} in which we provide a more explicit representation of the quadratic form for the effective operator in two dimensions.
\section{Preliminaries}\label{sec:prelimi}
In \S \ref{subsec:eqnDiracmat} we construct the Dirac matrices $(\alpha_j)_{j\in\{1,\dots,n+1\}}$. \S \ref{subsec:diffgeom} gathers various results of differential geometry that are used later on in the paper. \S \ref{subsec:transdira} is about the spectral properties of a transverse Dirac operator which plays an important role in the following.
\subsection{Representation of Clifford algebras}\label{subsec:eqnDiracmat}
The objective of this paragraph is, for a given $n\in \N$, to build the explicit hermitian matrices $\alpha_j \in \CC^{N\times N}$, where $N = 2^{\lfloor\frac{n+1}2\rfloor}$, $j \in \{1,\dots,n\}$ mentioned in \eqref{eqn:anticomalpha}. This is the purpose of the following proposition.
\begin{prop}\label{prop:alpha_Gamma}
	 Let $n\in \N$ and set $N = 2^{\lfloor\frac{n+1}2\rfloor}$. There exists a family $(\alpha_j)_{j\in \{1,\dots,n+1\}}$ of hermitian matrices of size $N\times N$, which satisfy the anti-commutation relations:
\begin{equation}
	\alpha_j \alpha_k + \alpha_k \alpha_j = 2 \delta_{j,k} I_{N},\quad \text{for } j,k \in\{1,\dots,n\}.
	\label{eqn:cacrel}
\end{equation}
Moreover, if for $x \in \RR^n$ we define $\Gamma(x) := \sum_{j=1}^n x_j \alpha_j$, for $x,y\in\dR^n$ there holds
\begin{equation}
	\Gamma(x)\Gamma(y) + \Gamma(y)\Gamma(x) = 2 (x\cdot y) I_N
	\label{eqn:Gammacomm}
\end{equation}
and the family $(\alpha_j)_{j\in\{1,\dots,n+1\}}$ can be chosen in such a way that $\Gamma(x)$ writes
\[
	\Gamma(x) = \begin{pmatrix}\mathbf{0}_{\frac{N}2}& \beta(x)^*\\\beta(x) & \mathbf{0}_{\frac{N}2}\end{pmatrix},
\]
where the maps $x\mapsto\beta(x)$ and $x\mapsto \beta(x)^*$ are linear and verify for all $x,y \in \RR^n$
\begin{equation}
	\beta(x)^*\beta(y) + \beta(y)^*\beta(x) = \beta(x)\beta(y)^*+ \beta(y)\beta(x)^* = 2 (x\cdot y)I_{\frac{N}2}.
	\label{eqn:betacomm}
\end{equation}
\label{prop:cla-main-prop}
\end{prop}
We start by proving the following Lemma, reminiscent of  \cite[Lemma 2.4]{MOP20} and following \cite[Chapter 15]{DG13} or \cite[Appendix E]{dWS86}.
\begin{lem} Let $n\in \N$. There exists a family of hermitian matrices $(\gamma_j(n))_{j\in\{1,\dots,n\}}$ of size $2^{\lfloor\frac{n}2\rfloor}\times 2^{\lfloor\frac{n}2\rfloor}$ such that for all $j,k \in \{1,\dots,n\}$ there holds
\[
	\gamma_j(n)\gamma_k(n) + \gamma_k(n)\gamma_j(n) = 2 \delta_{j,k} I_{2^{\lfloor\frac{n}2\rfloor}}.
\]
\label{lem:cliff}
\end{lem}
\begin{proof}
The argument proceeds by constructing recursively (on the dimension) these matrices.	
For $n = 1$, we set $\gamma_1(1) = 1$ and for $n=2$, we set
\[
	\gamma_1(2) = \begin{pmatrix}0 & 1\\1 & 0\end{pmatrix},\quad \gamma_1(2) = \begin{pmatrix} 0 & -i\\i  &0\end{pmatrix},
\]
where one recognizes the standard Pauli matrices. In particular, they are hermitian and for $j,k \in \{1,2\}$, one can check that there holds
\[
	\gamma_j(2)\gamma_k(2) + \gamma_k(2)\gamma_j(2) = 2 \delta_{j,k} I_2.
\]
Let $n \geq 3$ and assume we have constructed a family $(\gamma_j(n-1))_{j\in\{1,\dots,n-1\}}$ of hermitian matrices of size $2^{\lfloor\frac{n-1}2\rfloor}\times2^{\lfloor\frac{n-1}2\rfloor}$ such that for all $j,k \in \{1,\dots,n-1\}$ the following anti-commutation relation holds
\begin{equation}
	\gamma_j(n-1)\gamma_k(n-1) + \gamma_k(n-1)\gamma_j(n-1) = 2 \delta_{j,k} I_{2^{\lfloor\frac{n-1}2\rfloor}}.
	\label{eqn:antic-HR}
\end{equation}
Now, we distinguish two cases depending on whether $n$ is even or $n$ is odd.\\

\paragraph{\bf If $n = 2d+2$ (with $d\in \N$).} We set for all $j \in \{1,2,\dots,2d+1\}$
	\begin{equation}
		\gamma_j(2d+2) = \begin{pmatrix}\mathbf{0}_{2^d}& \gamma_j(2d+1)\\ \gamma_j(2d+1) & \mathbf{0}_{2^d}\end{pmatrix}
		\label{eqn:cons1}
	\end{equation}
and
\[
	\gamma_{2d+2}(2d+2) = \begin{pmatrix}\mathbf{0}_{2^d}& -i I_{2^d}\\ i I_{2^d} & \mathbf{0}_{2^d}\end{pmatrix}.
\]
For all $j \in \{1,\dots,2d+2\}$ the matrix $\gamma_j(2d+2)$ is hermitian and using the anti-commutation relations \eqref{eqn:antic-HR} satisfied by the matrices $(\gamma_j(2d+1))_{j\in\{1,\dots,2d+1\}}$ and the specific structure of the matrices $(\gamma_j(2d+2))_{j\in\{1,\dots,2d+2\}}$ one gets for all $j,k\in\{1,\dots,2d+2\}$ the anti-commutation relation
\[
	\gamma_j(2d+2)\gamma_k(2d+2) + \gamma_k(2d+2)\gamma_j(2d+2) = 2 \delta_{j,k} I_{2^{d+1}}.
\]
\paragraph{\bf If $n = 2d+1$ (with $d\in \N$).} We set for all $j \in \{1,\dots,2d\}$
\[
	\gamma_j(2d+1) = \gamma_j(2d)
\]
and
\[
	\gamma_{2d+1}(2d+1) = \begin{pmatrix}I_{2^{d-1}}&\mathbf{0}_{2^{d-1}} \\ \mathbf{0}_{2^{d-1}} & -I_{2^{d-1}} \end{pmatrix}.
\]
For all $j \in \{1,\dots,2d+1\}$ the matrix $\gamma_j(2d+1)$ is hermitian and using the anti-commutation relations \eqref{eqn:antic-HR} satisfied by the matrices $(\gamma_j(2d))_{j\in\{1,\dots,2d\}}$ one gets for all $j,k\in\{1,\dots,2d\}$ the anti-commutation relation
\begin{equation}
	\gamma_j(2d+1)\gamma_k(2d+1) + \gamma_k(2d+1)\gamma_j(2d+1) = 2 \delta_{j,k} I_{2^{d}}.
	\label{eqn:cliffrec-HRbis}
\end{equation}
One gets readily that $\gamma_{2d+1}^{2}(2d+1) = I_{2^d}$ and the only thing left to prove is that for all $j\in\{1,\dots,2d\}$ there holds
\[
	\gamma_j(2d+1)\gamma_{2d+1}(2d+1) + \gamma_{2d+1}(2d+1)\gamma_j(2d+1) = \mathbf{0}_{2^d}.
\]
Note that if $j \in\{1,\dots,2d-1\}$ there holds
\begin{align*}
	\gamma_j(2d+1)\gamma_{2d+1}(2d+1) &= \gamma_j(2d)\gamma_{2d+1}(2d+1)\\&= \begin{pmatrix} \mathbf{0}_{2^{d-1}} & \gamma_j(2d-1)\\ \gamma_j(2d-1) & \mathbf{0}_{2^{d-1}}\end{pmatrix}\begin{pmatrix}I_{2^{d-1}}&\mathbf{0}_{2^{d-1}} \\ \mathbf{0}_{2^{d-1}} & -I_{2^{d-1}} \end{pmatrix}\\
	& = \begin{pmatrix}\mathbf{0}_{2^{d-1}} & -\gamma_j(2d-1)\\ \gamma_j(2d-1) & \mathbf{0}_{2^{d-1}} \end{pmatrix}\\
	& = - \begin{pmatrix}I_{2^{d-1}}&\mathbf{0}_{2^{d-1}} \\ \mathbf{0}_{2^{d-1}} & -I_{2^{d-1}} \end{pmatrix}\begin{pmatrix} \mathbf{0}_{2^{d-1}} & \gamma_j(2d-1)\\ \gamma_j(2d-1) & \mathbf{0}_{2^{d-1}}\end{pmatrix}\\
	& = - \gamma_{2d+1}(2d+1)\gamma_{j}(2d+1)
\end{align*}
and we have
\begin{align*}
	\gamma_{2d}(2d+1)\gamma_{2d+1}(2d+1) &= \gamma_{2d}(2d) \gamma_{2d+1}(2d+1)\\
	& = \begin{pmatrix}\mathbf{0}_{2^{d-1}} & -i I_{2^{d-1}}\\i I_{2^{d-1}} & \mathbf{0}_{2^{d-1}}\end{pmatrix}\begin{pmatrix}I_{2^{d-1}}&\mathbf{0}_{2^{d-1}} \\ \mathbf{0}_{2^{d-1}} & -I_{2^{d-1}} \end{pmatrix}\\
	& = \begin{pmatrix}\mathbf{0}_{2^{d-1}} & i I_{2^{d-1}}\\i I_{2^{d-1}} & \mathbf{0}_{2^{d-1}}\end{pmatrix}\\
	& = - \begin{pmatrix}I_{2^{d-1}}&\mathbf{0}_{2^{d-1}} \\ \mathbf{0}_{2^{d-1}} & -I_{2^{d-1}} \end{pmatrix}\begin{pmatrix}\mathbf{0}_{2^{d-1}} & -i I_{2^{d-1}}\\i I_{2^{d-1}} & \mathbf{0}_{2^{d-1}}\end{pmatrix}\\
	& = - \gamma_{2d+1}(2d+1)\gamma_{2d}(2d+1).
\end{align*}
Thus, we have proved that Equation \eqref{eqn:cliffrec-HRbis} extends to all $j,k \in\{1,\dots,2d+1\}$, that is for all $j,k \in\{1,\dots,2d+1\}$ there holds
\[
	\gamma_j(2d+1)\gamma_k(2d+1) + \gamma_k(2d+1)\gamma_j(2d+1) = 2 \delta_{j,k} I_{2^{d}}.
\]
This achieves the construction by recursion.
\end{proof}

Now, we have all the tools to prove Proposition \ref{prop:cla-main-prop}.
\begin{proof}[Proof of Proposition \ref{prop:cla-main-prop}]
For all $j\in\{1,\dots,n+1\}$ we define the Dirac matrix $\alpha_j$ as
\begin{equation}
	\alpha_j := \left\{\begin{array}{lcl}
					\gamma_j(n+1) & \text{if} & n\text{ is even}\\
					\gamma_j(n+1) &\text{if} & n\text{ is odd and } j\in\{1,\dots,n\}\\
					\begin{pmatrix}I_{2^{\frac{n-1}2}} & \mathbf{0}_{2^{\frac{n-1}2}}\\
					\mathbf{0}_{2^{\frac{n-1}2}} & - I_{2^{\frac{n-1}2}}\end{pmatrix}&\text{if} & n\text{ is odd and } j=n+1
				\end{array}\right.,
			\label{eqn:defalpha}
\end{equation}
where the matrices $(\gamma_j(n+1))_{j\in\{1,\dots,n+1\}}$ are introduced in Lemma \ref{prop:cla-main-prop}.
By construction, when $n$ is even the family $(\alpha_{j})_{j\in\{1,\dots,n+1\}}$ is a family of hermitian matrices of size $N\times N$ satisfying the anti-commutation relations \eqref{eqn:cacrel}. Similarly, when $n$ is odd the family $(\alpha_{j})_{j\in\{1,\dots,n+1\}}$ is a family of hermitian matrices of size $N$ and for all $j,k \in \{1,\dots,n\}$ we have
\[
	\alpha_j\alpha_k + \alpha_k\alpha_j = 2\delta_{j,k}I_{N}.
\]
The only thing left to prove is that this relation also holds for $j,k \in\{1,\dots,n+1\}$. This is true because an elementary computation yields $\alpha_{n+1}^2 = I_N$ and for any $j\in \{1,\dots,n\}$ the matrix $\alpha_j = \gamma_j(n+1)$ is of the shape
\[
	\gamma_j(n+1) = \begin{pmatrix}\mathbf{0}_{2^{\frac{n-1}2}}& A\\B & \mathbf{0}_{2^{\frac{n-1}2}}\end{pmatrix},
\]
for some complex matrices $A,B$ of size $2^{\frac{n-1}2}\times2^{\frac{n-1}2}$ (see \eqref{eqn:cons1}). There holds
\[
	\begin{pmatrix}\mathbf{0}_{2^{\frac{n-1}2}}& A\\B & \mathbf{0}_{2^{\frac{n-1}2}}\end{pmatrix} \alpha_{n+1} = \begin{pmatrix}\mathbf{0}_{2^{\frac{n-1}2}}& -A\\B & \mathbf{0}_{2^{\frac{n-1}2}}\end{pmatrix}  = -\alpha_{n+1}\begin{pmatrix}\mathbf{0}_{2^{\frac{n-1}2}}& A\\B & \mathbf{0}_{2^{\frac{n-1}2}}\end{pmatrix}
\]
and the anti-commutation relations \eqref{eqn:cacrel} are satisfied by the family of Dirac matrices $(\alpha_j)_{j\in\{1,\dots,n+1\}}$.

Now, pick $x,y \in \RR^n$, one remarks that there holds
\begin{align*}
	\Gamma(x) \Gamma(y) = \sum_{j,k = 1}^n x_j y_k \alpha_j \alpha_k &= (x \cdot y) I_N + \sum_{j,k = 1\atop j\neq k}^n x_jy_k \alpha_j\alpha_k\\
	& = (x \cdot y) I_N - \sum_{j,k = 1\atop j\neq k}^n x_jy_k \alpha_k\alpha_j\\
	& = 2(x\cdot y) I_N - \sum_{j,k = 1}^n x_jy_k \alpha_k\alpha_j\\
	& = 2(x\cdot y) I_N - \Gamma(y) \Gamma(x),
\end{align*}
where we have used the anti-commutation relations \eqref{eqn:cacrel} and it yields \eqref{eqn:Gammacomm}.

Remark that with definition \eqref{eqn:defalpha}, for all $x \in \RR^n$ there holds
\[
	\Gamma(x) = \left\{\begin{array}{lcl}\begin{pmatrix}\mathbf{0}_{2^{d}} & \upsilon(x),\\\upsilon(x) & \mathbf{0}_{2^{d}}\end{pmatrix} & \text{if} & n = 2d+1\\\begin{pmatrix}\mathbf{0}_{2^{d-1}} & \lambda(x)^*\\\lambda(x) & \mathbf{0}_{2^{d-1}}\end{pmatrix} & \text{if} & n = 2d.\end{array}\right.,
\]
with
\[
	\upsilon(x) := \sum_{j=1}^{2d+1}x_j \gamma_j(2d+1),\qquad\lambda(x) := \sum_{j=1}^{2d-1}x_j \gamma_j(2d-1) + i x_{2d} I_{2^{d-1}}.
\]
and if one sets
\[
	\beta(x) := \left\{\begin{array}{lcl}
					\upsilon(x) & \text{if} & n \text{ is odd},\\
					\lambda(x) & \text{if} & n \text{ is even},
				\end{array}\right.
\]
we remark that there holds
\[
	\Gamma(x) = \begin{pmatrix} \mathbf{0}_{\frac{N}2} & \beta(x)^*\\ \beta(x) & \mathbf{0}_{\frac{N}2}\end{pmatrix}.
\]
Hence, for $x,y \in \RR^n$, \eqref{eqn:Gammacomm} gives

\begin{equation}
	\beta(x)^*\beta(y) + \beta(y)^*\beta(x) = \beta(x)\beta(y)^* + \beta(y)\beta(x)^* = 2(x\cdot y)I_{\frac{N}2}
\end{equation}
which is precisely \eqref{eqn:betacomm} and Proposition \ref{prop:cla-main-prop} is proved.

\end{proof}

\subsection{Standard tools from differential geometry}\label{subsec:diffgeom}

In this paragraph we collect several well known facts of differential geometry. The first lemma concerns the metric tensor $G$ associated with the map \eqref{eqn:tubcoord}.

\begin{lem} Let $(e_1,\dots,e_{n-1})$ be a local coordinate system at $s \in \Sigma$. If the shape operator $S$ is identified with its matrix in the basis $(e_1,\dots,e_{n-1})$ then, in the basis $(e_1,\dots,e_{n-1},\nu)$ the metric tensor $G$ writes
\begin{equation}
	G = \begin{pmatrix}g (I_{n-1} + \eps t S)^2 & 0 \\ 0 & \eps^2\end{pmatrix}
	\label{eqn:blockstructmetric}
\end{equation}
where, with a slight abuse of notation, we identify the metric $G$ with its matrix in the basis $(e_1,\dots,e_{n-1})$. Here $g = (g_{i,j})_{i,j\in\{1,\dots,n-1\}}$ is the metric tensor of the metric on $\Sigma$ induced by its embedding into $\RR^n$ expressed in the local coordinate system $(e_1,\dots,e_{n-1})$. In particular, the following holds.
\begin{enumerate}
\item\label{itm:1-metric} We have $\det(G) = \eps^2\det(g) \Big(1 +\sum_{k=1}^{n-1}\eps^k t^k H_k\Big)^2$.
\item\label{itm:3-metric} There exist $\varepsilon_0 > 0$ and  $c > 0$ such that for all $\varepsilon \in (0,\varepsilon_0)$ and $x = (x_1,\dots,x_{n-1},0)^\top \in \RR^n$ there holds
\[
	(1-c\eps) (g^{-1} x') \cdot x' \leq (G^{-1} x)\cdot x \leq (1+c\eps) (g^{-1} x') \cdot x',
\]
where $x' = (x_1,\dots,x_{n-1})^{\top}\in \RR^{n-1}$.
\end{enumerate}
\label{lemma:metric}
\end{lem}
\begin{proof} For $s\in\Sigma$, let $(e_1,\dots,e_{n-1})$ denote a local coordinate system at $s \in \Sigma$. By definition, for all $(s,t)\in \Sigma\times(-1,1)$ and all $j\in\{1,\dots,n-1\}$ there holds
\begin{equation}
	\partial_j \phi_\eps (s,t) = e_j + \eps t \partial_j \nu = e_j + \eps t \nabla_{e_j} \nu =(I_{n-1} + \eps t S) e_j \in T_s \Sigma
	\label{eqn:meblock1}
\end{equation}
and
\begin{equation}
	\partial_t \phi_\eps(s,t) = \eps \nu.
	\label{eqn:meblock2}
\end{equation}
Equations \eqref{eqn:meblock1} and $\eqref{eqn:meblock2}$ yield the expected block structure for the metric tensor $G$. Moreover, for all $j,k\in\{1,\dots,n-1\}$, there holds
\begin{align*}
	G_{j,k} = \partial_j \phi_\eps \cdot \partial_k \phi_\eps &= \big((I_{n-1}+ \eps t S)e_j \big)\cdot \big((I_{n-1}+ \eps t S)e_k\big)\\ &= e_j \cdot \big( (I_{n-1}+ \eps t S)^2e_k\big)\\ & =  e_j \cdot \Big( \sum_{p=1}^{n-1}((I_{n-1}+ \eps t S)^2)_{k,p}e_p\Big)\\&= \sum_{p=1}^{n-1}(e_j \cdot e_p)((I_{n-1}+ \eps t S)^2)_{k,p}\\&=\sum_{p=1}^{n-1}(e_j \cdot e_p)((I_{n-1}+ \eps t S)^2)_{p,k}\\&= \sum_{p=1}^{n-1}g_{j,p}	((I_{n-1}+ \eps t S)^2)_{p,k}\\& = \big(g(I_{n-1}+ \eps t S)^2)\big)_{j,k},
\end{align*}
where we have used twice the symmetry of $S$. It gives the block structure of $G$ stated in Equation \eqref{eqn:blockstructmetric}.

Regarding Point \eqref{itm:1-metric}, \eqref{eqn:blockstructmetric} gives
\[
	\det(G) = \eps^2 \det(g) \big(\det(I_{n-1} + \eps t S)\big)^2 = \eps^2\det(g)\Big(\prod_{k=1}^{n-1}(1+\eps t \kappa_k)\Big)^2,
\]
where we have used that $S$ is diagonalizable with the principal curvatures $(\kappa_k)_{k\in\{1,\dots,n-1\}}$ as eigenvalues. The only thing left to note is that
\[
	\prod_{k=1}^{n-1}(1+\eps t \kappa_k) = 1 + \sum_{k=1}^{n-1}\eps^k t^k H_k,
\]
which proves Point \eqref{itm:1-metric}.

Let us move to the proof of Point \eqref{itm:3-metric}. It is enough to prove it for $(e_1,\dots,e_{n-1})$ an orthonormal
local coordinate system at $s\in \Sigma$ which diagonalizes the shape operator $S$. In this basis the metric tensor $G$ is diagonal and writes
\[
	G_{j,j} = (1+\eps t \kappa_j)^{2}, \qquad \forall j \in \{1,\dots,n-1\}.
\] 

Now, pick $x = (x_1,\dots,x_{n-1},0)^\top \in \RR^n$ and set $x' = (x_1,\cdots,x_{n-1})^\top \in \RR^{n-1}$.
There holds:
\begin{align*}
	(G^{-1} x)\cdot x &= \sum_{j=1}^{n-1}\frac{1}{(1+\eps t \kappa_j)^2} |x_j|^2
\end{align*}
and we get
\begin{align*}
	\frac1{(1 + \eps \kappa_{max})^2} (g^{-1}x')\cdot x' = \frac1{(1 + \eps \kappa_{max})^2} |x'|^2 &\leq (G^{-1}x) \cdot x\\&\leq  \frac1{(1 - \eps \kappa_{max})^2}  |x'|^2 = \frac1{(1 - \eps \kappa_{max})^2}  (g^{-1}x')\cdot x',
\end{align*}
with $\kappa_{max} = \max_{s\in \Sigma}\{|\kappa_1(s)|,\dots,|\kappa_{n-1}(s)|\}$. This concludes the proof of Lemma \ref{lemma:metric}.
\end{proof}

The second lemma provides an expression of the mean curvature $H_1^{\partial\Omega_\eps}$ on $\partial \Omega_\eps$ in terms of the mean curvature $H_1$ of $\Sigma$. To state it, remark that there holds
\begin{equation}
	\partial\Omega_\eps := \Sigma_\eps^+ \cup \Sigma_\eps^-,
	\label{eqn:defsigmapm}
\end{equation}
where we have set $\Sigma_\eps^\pm := \{\phi_\eps(s,\pm1)\colon s \in \Sigma\}$ with $\phi_\eps$ being the map defined in \eqref{eqn:tubcoord}. The following lemma holds.

\begin{lem} For all $s \in \Sigma$, we have:
\[
	H_1^{\partial\Omega_\eps}(\phi_\eps(s,\pm1)) = \pm H_1(s) - \eps \big(H_1(s)^2 - 2H_2(s)\big)+ \mathcal{O}(\eps^2),
\]
where the remainder is understood in the $L^\infty(\Sigma)$-norm.
\label{lem:expcurv}
\end{lem}
\begin{proof}
	Let us fix $s \in \Sigma$. Let $\nu^\pm$ be the unit normals on $\Sg^\pm_\eps$ pointing outwards of $\Omg_\eps$. The principal curvatures of $\Sg^\pm_\eps$ at the point $s\pm\eps\nu(s)$ are given by the eigenvalues of the map in $T_{\phi_\eps(s,\pm1)}\Sg_\eps^\pm$ defined by
\[	
\cY\mapsto \nabla_\cY\nu^\pm(\phi_\eps(s,\pm1)).
\]
Let $I\subset\dR$ be an interval around the origin. Consider a smooth map $\varphi\colon I\arr\Sg$ with $\varphi(0) = s$ and $\varphi'(0)= \cX$. Let us define the map $\psi_\pm\colon I\arr\Sg_\eps^\pm$,
$\psi_\pm(t) := \varphi(t) \pm \eps\nu(\varphi(t))$. We have
$\psi_\pm(0) = \phi_\eps(s,\pm1)$ and $\psi'_\pm(0) =\varphi'(0) \pm \eps(\nb_{\varphi'(0)}\nu)(s) = \cX \pm \eps(\nb_{\cX}\nu)(s)$. If we set $\cX = e_j$ where $(e_1,\dots,e_{n-1})$ is a local orthonormal basis at $s\in\Sg$ which diagonalizes the map $\cY\mapsto\nb_{\cY}\nu$ we get that
\[
	\psi'_\pm(0) = (1\pm\eps\kp_j)e_j =:f_j
\]
Then we find that
\begin{equation}\label{eq:der_fj_1}
\begin{aligned}
 \nb_{f_j}\nu^\pm(s\pm\eps\nu(s))& = \frac{\dd}{\dd t}\left[\nu^\pm(\psi_\pm(t))\right]\Big|_{t=0}\\ 
 &=
 \pm\frac{\dd}{\dd t}\left[\nu(\varphi(t))\right]\Big|_{t=0} = 
 \pm (\nb_{e_j}\nu)(s) = \pm\kp_j e_j,
\end{aligned}
\end{equation}
where we have used that for all $t \in I$ we have $\nu^\pm(\psi_\pm(t)) = \pm \nu(\varphi(t))$.
On the other hand, as the shape operator is a linear map, we have
\begin{equation}\label{eq:der_fj_2}
	\nb_{f_j}\nu^\pm(\phi_\eps(s,\pm1)) 
	= 
	(1\pm\eps\kp_j)\nb_{e_j}\nu^\pm(\phi_\eps(s,\pm1)).
\end{equation}
Combining~\eqref{eq:der_fj_1} and~\eqref{eq:der_fj_2} we obtain that
\[
	(\nb_{e_j}\nu^\pm)(\phi_\eps(s,\pm1)) = \pm\frac{\kp_j}{1\pm\eps\kp_j}e_j.
\]
Hence, we get the expressions for the mean curvatures of $\Sg_\eps^\pm$
\[
\begin{aligned}
H_1^{\partial\Omega_\eps}(\phi_\eps(s,\pm1)) 
	&= 
	\pm\sum_{j=1}^{n-1}
	\frac{\kp_j(s)}{1\pm\eps\kp_j(s)}\\
& = \pm H_1(s) -\eps\sum_{j=1}^{n-1}\kp_j^2 + \cO(\eps^2)\\
& = \pm H_1(s) - \eps(H_1^2(s) -2H_2(s))+\cO(\eps^2),
\end{aligned}
\]
which is precisely Lemma \ref{lem:expcurv}.
\end{proof}

\subsection{The transverse Dirac operator}\label{subsec:transdira}
Let $x \in \mathbb{S}^{n-1}\subset\dR^n$ be fixed.
Consider the following self-adjoint operator in the Hilbert space $L^2((-1,1),\CC^N)$:
\begin{equation}\label{eqn:deftransdirac}
\begin{aligned}
	\sfT_x f = \sfT_x(m) f &:= -i \Gamma(x) f' + m \alpha_{n+1}f,\\
	\dom{\sfT_x} &:= \{f \in H^1((-1,1),\CC^N) \colon -i \alpha_{n+1}\Gamma(x) f(\pm 1) =\pm f(\pm1)\},
\end{aligned}
\end{equation}
where the matrices $\Gamma(x)$ and $\alpha_{n+1}$ are introduced as in Proposition~\ref{prop:alpha_Gamma}.

The goal of this paragraph is to prove the following proposition.
\begin{prop}\label{prop:1doperator} The operator $\sfT_x$ is self-adjoint and has compact resolvent. Moreover, the following holds.
\begin{enumerate}
	\item \label{itm:1-1D}For all $f \in \dom{\sfT_x}$ there holds
	\[
		\|\sfT_x f\|_{L^2((-1,1),\CC^N)}^2 = \|f'\|_{L^2((-1,1),\CC^N)}^2 + m^2 \|f\|_{L^2((-1,1),\CC^N)}^2 + m \big(|f(1)|_{\CC^N}^2 + |f(-1)|_{\CC^N}^2\big).
	\]
	\item\label{itm:2-1D} For all $p \geq 1$, define $k_p(m)$ as the only root lying in $[(2p-1)\frac\pi4,p\frac\pi2]$ of $m\sin(2k) +k\cos(2k) = 0$. If one sets $E_p(m) := \sqrt{m^2 + k_p(m)^2}$ there holds
	\[
		Sp(\sfT_x) = \bigcup_{p\geq 1}\{ E_p(m),-E_p(m)\}.
	\]
	\item \label{itm:3-1D}When $m\to 0$, the following asymptotic expansion holds
	\[
		k_1(m) = \frac\pi4 + \frac2\pi m - \frac{16}{\pi^3}m^2  + \mathcal{O}(m^3).
	\]
	\item \label{itm:4-1D}For $p\geq 1$, $\dim(\ker(\sfT_x - E_p)) = \dim(\ker(\sfT_x + E_p)) = \frac{N}2$. If $(\varepsilon_j)_{j\in\{1,\dots,\frac{N}2\}}$ is the canonical basis of $\CC^{\frac{N}2}$, a normalized eigenfunction associated with $E_p(m)$ is given for all $t \in (-1,1)$ by
	\begin{multline}
		\varphi_{j,p}^{m,+}(t) := N_{m,p}\Bigg(k_p(m) \cos(k_p(m)(t+1))\begin{pmatrix}\varepsilon_j\\-i\beta(x)\varepsilon_j\end{pmatrix}\\+\sin(k_p(m)(t+1)) \begin{pmatrix}(E_p(m)+m)\varepsilon_j\\i(E_p(m)-m)\beta(x)\varepsilon_j\end{pmatrix}\Bigg)
		\label{eqn:posmod1D}
	\end{multline}
and a normalized eigenfunction associated with $-E_p(m)$ is given for all $t \in (-1,1)$ by
		\begin{multline}
		\varphi_{j,p}^{m,-}(t) := N_{m,p}\Bigg(k_p(m) \cos(k_p(m)(t+1))\begin{pmatrix}i\beta(x)^*\varepsilon_j\\\varepsilon_j\end{pmatrix}\\+\sin(k_p(m)(t+1)) \begin{pmatrix}i(-E_p(m)+m)\beta(x)^*\varepsilon_j\\(E_p(m)+m)\varepsilon_j\end{pmatrix}\Bigg),
		\label{eqn:negmod1D}
	\end{multline}
	where $N_{m,p}$ is a normalization constant. Moreover, the family $(\varphi_{j,p}^{m,\pm})_{j \in \{1,\dots,\frac{N}2\}}$ is an orthonormal basis of $\ker(\sfT_x \mp E_p(m))$.
	\item \label{itm:5-unitary} For any $x, y\in\dS^{n-1}$,
	$\sfT_y^2 = \Theta_{x,y}\sfT_x^2\Theta_{x,y}^*$
	with the unitary operator $\Theta_{x,y}$ in $L^2((-1,1),\dC^N)$ defined by \[
	(\Theta_{x,y} f)(t) := \frac12\big(I + i\Gamma(y)\big)\big(I-i\Gamma(x))\big)f(t),\qquad t\in(-1,1). 
	\] 
\end{enumerate}
\label{prop:1D}
\end{prop}
Before going through the proof of Proposition \ref{prop:1D} it is of crucial importance for further uses to note that for all $j \in \{1,\dots,\frac{N}2\}$ and $t \in (-1,1)$ there holds
	\begin{equation}
		\varphi_{j}^{+}(t) := \varphi_{j,1}^{0,+}(t) = \frac12 \cos\left(\frac\pi4(t+1)\right)\begin{pmatrix}\varepsilon_j\\-i\beta(x)\varepsilon_j\end{pmatrix}+\frac12\sin\left(\frac\pi4(t+1)\right) \begin{pmatrix}\varepsilon_j\\i\beta(x)\varepsilon_j\end{pmatrix}
		\label{eqn:modpos1D0}
	\end{equation}
and
	\begin{equation}
		\varphi_{j}^{-}(t) := \varphi_{j,1}^{0,-}(t) = \frac12 \cos\left(\frac\pi4(t+1)\right)\begin{pmatrix}i\beta(x)^*\varepsilon_j\\\varepsilon_j\end{pmatrix}+\frac12\sin\left(\frac\pi4(t+1)\right) \begin{pmatrix}-i\beta(x)^*\varepsilon_j\\\varepsilon_j\end{pmatrix}.
		\label{eqn:modneg1D0}
	\end{equation}
\begin{proof}[Proof of Proposition \ref{prop:1D}] Let us start by proving that $\sfT_x$ is a self-adjoint operator. As the multiplication operator by $\alpha_{n+1}$ is self-adjoint, it is enough to do so for $m =0$. Consider $f,g \in \dom{\sfT_x}$ and remark that
\begin{align*}
	\langle \sfT_x f, g\rangle_{L^2((-1,1),\CC^N)} &= \int_{-1}^1 \langle-i \Gamma(x) f'(t),g(t)\rangle_{\CC^N} \dd t \\&= \int_{-1}^1 \langle f(t),-i\Gamma(x)g'(t)\rangle_{\CC^N} \dd t \\& \qquad\qquad- i \big(\langle f(1),\Gamma(x) g(1)\rangle_{\CC^N} - \langle f(-1),\Gamma(x) g(-1)\rangle_{\CC^N}\big).
\end{align*}
Using the boundary conditions we obtain
\begin{align*}
	\langle f(1),\Gamma(x) g(1)\rangle_{\CC^N} &= \langle -i\alpha_{n+1}f(1),-i\alpha_{n+1}\Gamma(x) g(1)\rangle_{\CC^N}\\ &= \langle-i\alpha_{n+1}(-i\alpha_{n+1}\Gamma(x))f(1),g(1)\rangle_{\CC^N}\\& = - \langle\Gamma(x) f(1),g(1)\rangle_{\CC^N} = - \langle f(1),\Gamma(x) g(1)\rangle_{\CC^N},
\end{align*}
because $\Gamma(x)^* = \Gamma(x)$.
Hence, we get $\langle f(1),\Gamma(x) g(1)\rangle_{\CC^N} = 0$ and in the same way, one can prove that ${\langle f(-1),\Gamma(x) g(-1)\rangle_{\CC^N}= 0}$. Thus, $\sfT_x$ is a symmetric operator. Recall that
\begin{multline*}
	\dom{\sfT_x^*} := \big\{f \in L^2((-1,1),\CC^N) : \exists g \in L^2((-1,1),\CC^N) \text{ such that}\\ \forall h \in \dom{\sfT_x},\langle f, \sfT_x h\rangle_{L^2((-1,1),\CC^N)} = \langle g, h\rangle_{L^2((-1,1),\CC^N)}\big\}
\end{multline*}
and as for $f \in \dom \sfT_x^*$, there exists $g \in L^2((-1,1),\CC^N)$ such that $\langle f, \sfT_x h\rangle_{L^2((-1,1),\CC^N)} = \langle g, h\rangle_{L^2((-1,1),\CC^N)}$, we have by definition $\sfT_x^* f = g$. Let us show that $\dom{\sfT_x^*} \subset \dom{\sfT_x}$. Pick $f \in \dom{\sfT_x^*}$. There exists $g \in L^2((-1,1),\CC^N)$ such that for any $h \in \dom{\sfT_x}$ there holds
\[
	\langle f, \sfT_x h\rangle_{L^2((-1,1),\CC^N)} = \langle g, h\rangle_{L^2((-1,1),\CC^N)}.
\]
Choosing $h \in C_0^\infty((-1,1),\CC^N) \subset \dom\sfT_x$ and rewriting this equality in the sense of distributions yields $ g = - i \Gamma(x) f' \in L^2((-1,1),\CC^N)$ and in particular $f \in H^1((-1,1),\CC^N)$ and $\sfT_x^* f = - i \Gamma(x) f' $. Moreover, for all $h \in \dom{\sfT_x}$, an integration by parts yields
\begin{align*}
	\langle \sfT_x^* f, h\rangle_{L^2((-1,1),\CC^N)} &= \langle f, \sfT_x h\rangle_{L^2((-1,1),\CC^N)}\\&\qquad\qquad- i \big(\langle f(1),\Gamma(x) h(1)\rangle_{\CC^N} - \langle f(-1),\Gamma(x) h(-1)\rangle_{\CC^N}\big)
\end{align*}
which gives
\[
	\langle f(1),\Gamma(x) h(1)\rangle_{\CC^N} - \langle f(-1),\Gamma(x) h(-1)\rangle_{\CC^N} = 0.
\]
As this is true for all $h \in \dom{\sfT_x}$ we obtain that
\[
	\Gamma(x)f(1) \in \ker(-i\alpha_{n+1}\Gamma(x) - I_{N})^\perp \text{ and }\Gamma(x)f(-1) \in \ker(-i\alpha_{n+1}\Gamma(x) + I_{N})^\perp.
\]
Now, one remarks that $\ker(-i\alpha_{n+1}\Gamma(x) - I_{N})^\perp = \ker(-i\alpha_{n+1}\Gamma(x) + I_N)$ and $\ker(-i\alpha_{n+1}\Gamma(x) + I_{N})^\perp = \ker(-i\alpha_{n+1}\Gamma(x) - I_N)$. Indeed, we have the following orthogonal sum decomposition holds
\[
	\CC^N = \ker(-i\alpha_{n+1}\Gamma(x) + I_{N}) \oplus\ker(-i\alpha_{n+1}\Gamma(x) - I_{N})
\]
as can be seen decomposing any vector $v \in \CC^N$ as
\[
	v = \frac12(i\alpha_{n+1}\Gamma(x) + I_{N})v + \frac12(-i\alpha_{n+1}\Gamma(x) + I_{N})v := v_1 + v_2.
\]
and remarking that the anti-commutation relation $\Gamma(x)\alpha_{n+1} = - \alpha_{n+1}\Gamma(x)$ gives \\${v_1 \in \ker(-i\alpha_{n+1}\Gamma(x) + I_{N})}$, $v_2 \in \ker(-i\alpha_{n+1}\Gamma(x) - I_{N})$ and $\langle v_1,v_2\rangle_{\CC^N} = 0$.
Thus, we obtain
\[
	0 = (-i\alpha_{n+1}\Gamma(x) + I_N)\Gamma(x) f(1) = - i\alpha_{n+1} f(1) + \Gamma(x) f(1)
\]
which reads $-i\alpha_{n+1}\Gamma(x) f(1) = f(1)$. Similarly, one can prove that $i \alpha_{n+1}\Gamma(x) f(-1) = f(-1)$ which yields that $\sfT_x$ is self-adjoint. As $\dom{\sfT_x} \subset H^1((-1,1),\CC^N)$, the compact embedding of $H^1((-1,1),\CC^N)$ into $L^2((-1,1),\CC^N)$ gives the compactness of the resolvent of $\sfT_x$.

Let us prove Point \eqref{itm:1-1D}. To this aim pick $f\in \dom{\sfT_x}$. There holds
\[
\begin{aligned}
	&\|-i\G(x)f' + m\aa_{n+1} f\|^2 &= \|f'\|^2 + m^2\|f\|^2 + 2m\Re\Big(\int_{-1}^1\langle-i\G(x)f',\aa_{n+1}f\rangle_{\dC^N}\dd t\Big),
\end{aligned}
\]
where we use the abbreviation $\|\cdot\|$ for the norm in $L^2((-1,1),\CC^N)$.
Using an integration by parts we get
\[
\begin{aligned}
	\int_{-1}^1\langle-i\G(x)f',\aa_{n+1}f\rangle_{\dC^N}\dd t &=-\int_{-1}^1\langle-i\G(x) f, \aa_{n+1} f'\rangle_{\dC^N}\dd t +
	\left[\langle-i\G(x)f,\aa_{n+1}f\rangle_{\dC^N}\right]_{t=-1}^{t=1}\\
	&= \int_{-1}^1\langle\G(x)f, -i\aa_{n+1} f'\rangle_{\dC^N}\dd t\\
	&\quad+\Big[\langle-i\aa_{n+1}\G(x)f(1),f(1)\rangle_{\dC^N} -
	\langle-i\aa_{n+1}\G(x)f(-1),f(-1)\rangle_{\dC^N}
	\Big]\\
	&=-\int_{-1}^1\langle\aa_{n+1}f,-i\G(x)f'\rangle_{\dC^N}\dd t + |f(1)|^2_{\dC^N} + |f(-1)|^2_{\dC^N}.
\end{aligned}
\]
Hence, we get that
\[
	2\Re\left(\int_{-1}^1\langle-i\G(x)f',
	\aa_{n+1}f\rangle_{\dC^N}\dd t\right)=|f(1)|^2_{\dC^N} + |f(-1)|^2_{\dC^N}
\]
and therefore
\[
	\|\sfT_x f\|^2 = \|f'\|^2 +m^2\|f\|^2 +m\left(|f(-1)|^2_{\CC^N}+|f(1)|^2_{\CC^N}\right).
\]
which is precisely Point \eqref{itm:1-1D}.

Let us move to the proof of Point \eqref{itm:2-1D}. First, remark that by Point \eqref{itm:1-1D} and the min-max principle (see \eqref{eqn:minmax}), there holds $\mu_1(\sfT_x^2) \geq m^2$, where $\mu_1(\sfT_x^2)$ is the lowest eigenvalue of $\sfT_x^2$ according to the notation introduced in \eqref{eqn:notminmax}. In particular $Sp(\sfT_x) \cap (-m,m) = \emptyset$. Actually, $Sp(\sfT_x) \cap [-m,m] = \emptyset$. Indeed, if $\pm m$ is an eigenvalue of $\sfT_x$ associated with an eigenfunction $f$, Point \eqref{itm:1-1D} yields
\[
	m^2 \|f\|_{L^2((-1,1)),\CC^N)}^2 = \|\sfT_x f\|_{L^2((-1,1)),\CC^N)}^2 \geq \|f'\|_{L^2((-1,1)),\CC^N)}^2 +m^2\|f\|_{L^2((-1,1),\CC^N)}^{2}
\]
thus $f' = 0$ and $f$ is a constant function. This is impossible because in order to satisfy the boundary conditions $f$ needs to be zero and we have $Sp(\sfT_x) \cap [-m,m] = \emptyset$. Now, let $E$ be an eigenvalue of $\sfT_x$ associated with $f = (f_1,f_2)^\top \in \dom{\sfT_x}$, where $f_1,f_2\in H^1((-1,1),\CC^{N/2})$. Hence, the eigenvalue equation reads
\begin{equation}
	\left\{
		\begin{array}{lcl}
			-i \beta(x)^* f_2' +mf_1 = E f_1\\
			-i\beta(x) f_1' - m f_2 = E f_2
		\end{array},
	\right.
	\label{eqn:sysev1D}
\end{equation}
in particular, $-f_1'' = k^2 f_1$, where we have set $k = \sqrt{E^2- m^2}$ and used that $\beta(x)^*\beta(x) = I_{\frac{N}2}$ (see Equation \eqref{eqn:betacomm}). Thus, for all $t \in (-1,1)$ there holds
\[
	f_1(t) = \cos(k(t+1)) A + \sin(k(t+1)) B
\]
for some vectors $A,B \in \CC^{\frac{N}2}$. The second equation in \eqref{eqn:sysev1D} gives
\[
	f_2(t) = \frac{-i k }{E+m}\cos(k(t+1))\beta(x) B + \frac{ik}{E+m}\sin(k(t+1))\beta(x) A
\]
which rewrites
\[
	f(t) = \cos(k(t+1)) \begin{pmatrix}A \\\frac{-i k }{E+m}\beta(x) B\end{pmatrix} + \sin(k(t+1)) \begin{pmatrix} B \\ \frac{ik}{E+m}\beta(x) A\end{pmatrix}.
\]
Now, let us deal with the boundary conditions. Remark that we have
\[
	-i\alpha_{n+1} \Gamma(x) = \begin{pmatrix} \mathbf{0}_{\frac{N}2} & -i \beta(x)^* \\ i \beta(x) &\mathbf{0}_{\frac{N}2}\end{pmatrix}.
\]
Hence, the boundary conditions read $f_2(\pm1) = \pm i \beta(x) f_1(\pm1)$. For $t = -1$ it yields
\[
	A = \frac{k}{E+m} B
\]
and $f$ rewrites for all $t\in(-1,1)$ as
\[
	f(t) = \frac{k}{E+m}\cos(k(t+1)) \begin{pmatrix} B \\-i\beta(x) B\end{pmatrix} + \sin(k(t+1)) \begin{pmatrix} B \\ i\frac{E-m}{E+m}\beta(x) B\end{pmatrix},
\]
where we have used that $k^2 = E^2-m^2 = (E-m)(E+m)$.
Choosing $B = (E+m) C$ for some $C \in \CC^{\frac{N}2}$, there holds
\begin{equation}
	\label{eqn:fp1D}
	f(t) = k \cos(k(t+1))\begin{pmatrix}C \\ -i\beta(x) C\end{pmatrix} + \sin(k(t+1))\begin{pmatrix}(E+m) C\\ i(E-m) \beta(x) C\end{pmatrix}.
\end{equation}
The boundary condition at $t = 1$ yields
\begin{equation}
	k \cos(2k) + m\sin(2k) = 0.
	\label{eqn:implicit}
\end{equation}
The solution $k = 0$ is excluded because otherwise $E = \pm m$ which is not possible. Hence, we are looking for positive solutions of Equation \eqref{eqn:implicit}. It can be checked that for $p\in \N$, this equation has a unique solution denoted $k_p(m)$ in the interval $[(2p-1)\frac\pi4,p\frac\pi2]$ which proves Point \eqref{itm:2-1D}.

To prove Point \eqref{itm:3-1D}, we use the implicit function theorem near the point $(\frac\pi4,0)$ for the function $G(k,m) := k\cos(2k) + m \sin(2k)$. Indeed, $G \in C^\infty(\RR^2)$ and there holds $G(\frac\pi4,0) = 0$ as well as $\partial_k G(\frac\pi4,0) = - \frac\pi2 \neq 0$ hence, there exists $\delta_1,\delta_2 > 0$ and $K : (-\delta_1,\delta_1) \to (\frac\pi4-\delta_2,\frac\pi4+\delta_2)$ such that for all $m \in (-\delta_1,\delta_1)$ there holds $G(K(m),m) = 0$, $K(0) = \frac\pi4$ and $K \in C^\infty((-\delta_1,\delta_1))$. One remarks that for $m \in [0,\delta_1)$ there holds $K(m) = k_1(m)$. In particular, when $m\to 0$, we have
\[
	k_1(m) = K(0) + K'(0)m + \frac12K''(0)m^2 + \mathcal{O}(m^3).
\]
Using the equation $G(K(m),m) = 0$ for all $m \in (-\delta_1,\delta_1)$ we get
\[
	K'(0) = \frac2\pi,\quad K''(0) = -\frac{32}{\pi^3}
\]
which proves Point \eqref{itm:3-1D}.

Point \eqref{itm:4-1D} is a consequence of \eqref{eqn:fp1D}. Indeed, any $C \in \CC^{\frac{N}2}$ gives rise to an eigenfunction and there is an isomorphism between $\CC^{N/2}$ and $\ker(\sfT_x \mp E_p)$ which yields $\dim(\ker(\sfT_x \mp E_p)) = \frac{N}{2}$. The expression for $\varphi_{j,p}^{m,+}$ is obtained by choosing $C = \varepsilon_j$ and normalizing $f$ in Equation \eqref{eqn:fp1D}. The expression for $\varphi_{j,p}^{m,-}$ is obtained by choosing $C = i \beta(x)^* \varepsilon_j$ and normalizing $f$ in Equation \eqref{eqn:fp1D}. The orthogonality of the family $(\varphi_{j,p}^{m,\pm})_{j\in \{1,\dots,\frac{N}2\}}$ comes from the orthogonality in $\CC^{\frac{N}2}$ of the family $(\varepsilon_j)_{j\in \{1,\dots,\frac{N}2\}}$, which proves Point \eqref{itm:4-1D}. 

In order to prove Point~\eqref{itm:5-unitary}, we fix $x,y \in \mathbb{S}^{n-1}$ and start by checking that the operator $\Theta_{x,y}$ is unitary which is equivalent to the fact that the matrix $U = \frac12(I+i\Gamma(x))(I-i\Gamma(y))$ is unitary in $\dC^N$. Indeed we find that
\[
\begin{aligned}
	U^*U 
	&= 
	\frac14\big(I+i\Gamma(y)\big)
	\big(I-i\Gamma(x)\big)
	\big(I+i\Gamma(x)\big)
	\big(I-i\Gamma(y)\big)\\
	&=
	\frac12\big(I+i\Gamma(y)\big)
	\big(I-i\Gamma(y)\big)= I.	
\end{aligned}\]
Analogously we check that $UU^* = I$.
Using the matrix identities
\[
\begin{aligned}
	U^*\alpha_{n+1}\G(x)U 
	&=
	\frac14
	\big(I+i\G(y)\big)\big(I-i\G(x)\big)
	\alpha_{n+1} \G(x)\big(I+i\G(x)\big)\big(I-i\G(y)\big)\\
	& = 
	\frac14
	\alpha_{n+1}
	\big(I-i\G(y)\big)\big(I+i\G(x)\big) \G(x)\big(I+i\G(x)\big)\big(I-i\G(y)\big)\\
	&=
	\frac14
	\alpha_{n+1}i
	\big(I-i\G(y)\big)\big(I-i\G(x)\big) \big(I+i\G(x)\big)\big(I-i\G(y)\big)\\
	&=\frac12\alpha_{n+1}i\big(I-i\G(y)\big)^2 = \alpha_{n+1}\G(y),
	\end{aligned}
\]
we arrive at
\[
\begin{aligned}
	\Theta_{x,y}(\dom \sfT_x)& = \big\{f\in H^1((-1,1),\dC^N)\colon -i\aa_{n+1}\G(x)Uf(\pm 1) = \pm U f(\pm 1)\big\}\\
	& =\big\{f\in H^1((-1,1),\dC^N)\colon -i\aa_{n+1}\G(y)f(\pm 1) = \pm f(\pm 1)\big\} = \dom\sfT_y,
\end{aligned}
\]
and furthermore for $f \in \dom \sfT_y$ we get
\[
	\begin{aligned}
	\|\sfT_x\Theta_{x,y}^*f\|^2_{L^2((-1,1),\dC^N)} &\!=\! \| U f'\|^2_{L^2((-1,1),\dC^N)} + m^2\|U f\|^2_{L^2((-1,1),\dC^N)} + m\big(|Uf(1)|_{\dC^N}^2 + |Uf(-1)|_{\dC^N}^2\big)\\
	&\!=\! \|\sfT_y f\|^2_{L^2((-1,1),\dC^N)}. 
	\end{aligned}
\]
Thus, we conclude that $\sfT_y^2 = \Th_{x,y}\sfT_x^2\Th_{x,y}^*$.
\end{proof}

\section{Quadratic forms}\label{sec:qf}
In order to prove Theorem \ref{thm:main}, we work with the quadratic form associated with $\scD_\eps^2$. The following lemma is the starting point of our analysis and can be found, {\it e.g.}, in \cite{MOP20}.

\begin{lem}\cite[Lem. 2.1]{MOP20}
	For all $u\in\dom\scD_\eps$ there holds
	\begin{equation}
		\|\scD_\eps u\|^2_{L^2(\Omg_\eps,\CC^N)} = 
		\|\nb u\|^2_{L^2(\Omg_\eps,\CC^N)} + m^2\|u\|^2_{L^2(\Omg_\eps,\CC^N)} +\int_{\p\Omg_\eps}\left(m+\frac{H_1^{\partial\Omega_\eps}}{2}\right)|u|^2 \dd s
		\label{eqn:qfdsq}
	\end{equation}
	where $H_1^{\partial\Omega_\eps}$ is the mean curvature of $\p\Omg_\eps$ and $\dd s$ is the $(n-1)$-dimensional Hausdorff measure on $\partial\Omega_\varepsilon$.
\end{lem}
The purpose of this section is to prove the following proposition.
\begin{prop}\label{prop:quadratic_forms} There exist $\eps_0 > 0$, $c > 0$ such that for all $\eps \in (0,\eps_0)$ there is a unitary map $\cU : L^2(\Omega_\eps,\CC^N) \to L^2(\Sigma\times(-1,1),\CC^N, \dd s \dd t)$  such that for all $u \in \dom{\scD_\eps}$ there holds
\[
	c_\eps^{-}[\cU u] \leq \|\scD_\eps u\|_{L^2(\Omega_\eps,\CC^N)}^2 \leq c_\eps^+[\cU u],
\]
where the quadratic forms $c_\eps^\pm$ are defined by
\[
\begin{aligned}
	c_\eps^\pm[w] &:= (1\pm c\eps)\int_{\Sg\tm(-1,1)}|\nb_\Sg w|^2\dd s \dd t + \int_{\Sg\tm(-1,1)} \left(H_2-\frac{H_1^2}{4}\right)|w|^2\dd s \dd t\\
&\quad +\frac{1}{\eps^2}\left[\int_{\Sg\tm(-1,1)}|\p_tw|^2 \dd s \dd t+ (m\eps \pm c\eps^3)\int_\Sg(|w(s,-1)|^2+|w(s,1)|^2)\dd s\right]\\
&\qquad\qquad + (m^2\pm c\eps)\int_{\Sg\tm(-1,1)}|w|^2 \dd s \dd t,\\
\dom{c_\eps^\pm} &:= \big\{w \in H^1(\Sigma\times(-1,1),\CC^N) : \ -i \alpha_{n+1}\Gamma\big(\nu\big) w(\cdot\,,\pm1) = \pm w(\cdot\,,\pm1)\big\};
\end{aligned}
\]
here $\nabla_\Sg$ is the surface gradient on $\Sg$.
\label{prop:lbubqf}
\end{prop}
The proof amounts to rewrite the quadratic form given in \eqref{eqn:qfdsq} in tubular coordinates using the map $\phi_\eps$ introduced in \eqref{eqn:tubcoord}. In this set of coordinates, we obtain a weighted  $L^2$-space on an $n$-dimensional manifold, the weight being given by a metric. Then, by an adequate unitary transform, we go to a {\it flat metric} on $\Sigma\times (-1,1)$. We pay the price of new potential terms appearing in the quadratic form involving the previous metric. We handle each of these terms by straightforward (though demanding) expansions in powers of $\eps$.\\

\begin{proof}[Proof of Proposition \ref{prop:lbubqf}]
The proof of this proposition is split into three steps. In the first one, we rewrite the quadratic form with the help of the map $\phi_\eps$ introduced in \eqref{eqn:tubcoord}. Then, the resulting quadratic form acts in a weighted $L^2$-space and in the second step, we use a unitary map to get rid of this weight. The last step deals with the asymptotic expansion of several terms appearing in the expression of this new quadratic form.

\paragraph{Step 1.} Let us define the following unitary map
\begin{equation*}\label{key}
	\sfU\colon L^2(\Omg_\eps,\dC^N)\!\arr\! L^2(\Sg\times(-1,1),\dC^N; \varphi_\eps(s,t) \dd s \dd t),\qquad
	(\sfU u)(s,t)\!  :=\! u(\phi_\eps(s,t)),
\end{equation*}
where, if $G$ is the metric associated with the map $\phi_\eps$ defined in \eqref{eqn:tubcoord}, we have used that $\det(G) = \det(g) \varphi_\eps(s,t)^2$ (see Point \eqref{itm:1-metric}~Lemma \ref{lemma:metric}). Here for $(s,t)\in\Sigma\times[-1,1]$ we have set $\varphi_\eps(s,t) : = \eps\big(1+\sum_{k=1}^{n-1}\eps^k t^k H_k(s)\big)$ (it is again a consequence of Point \eqref{itm:1-metric}~Lemma \ref{lemma:metric}).
The quadratic form for the operator $\scD_\eps^2$ is given by
\[
	a_\eps[u] := \|\scD_\eps u\|^2_{L^2(\Omg_{\eps},\CC^N)},\qquad u \in \dom a_\eps := \dom\scD_\eps.
\]
We define the unitarily equivalent form in the Hilbert space $L^2(\Sg\times(-1,1),\dC^N; \varphi_\eps \dd s \dd t)$.
\[
	b_\eps[v] := a_\eps[\sfU^{-1} v],\qquad v \in \dom b_\eps := \sfU(\dom a_\eps).
\]
Using Equation \eqref{eqn:blockstructmetric} in Lemma \ref{lemma:metric}, for $v \in \dom b_\eps$, we find 
\[
\begin{aligned}
	\|\nb(\sfU^{-1} v)\|^2_{L^2(\Omg_\eps,\dC^N)} &= \int_{-1}^1\int_\Sg \sum_{j,k=1}^n G^{j,k}\langle \p_j v,\p_k v\rangle_{\CC^N} \varphi_\eps(s,t)\dd s \dd t\\
	&=\int_{-1}^1\int_\Sg\left\{\left(\sum_{j,k=1}^{n-1} G^{j,k}\langle\p_j v,\p_k v\rangle_{\CC^N}\right) +\frac{1}{\eps^2}|\p_t v|^2\right\} \varphi_\eps(s,t) \dd s \dd t
\end{aligned}\]
and using Point \eqref{itm:3-metric}~Lemma \ref{lemma:metric} we obtain that there exists $c > 0$ such that for all $\eps$ small enough there holds
\begin{multline}
	\int_{\Sigma\times(-1,1)}\Big((1-c\eps)|\nabla_\Sigma v|^2 + \frac1{\eps^2}|\partial_t v|^2  \Big)\varphi_\eps(s,t)\dd s \dd t\leq \|\nb(\sfU^{-1} v)\|^2_{L^2(\Omg_{\eps},\CC^N)} \\
	\leq \int_{\Sigma\times(-1,1)}\Big((1+c\eps)|\nabla_\Sigma v|^2 + \frac1{\eps^2}|\partial_t v|^2 \Big)\varphi_\eps(s,t) \dd s \dd t.
	\label{eqn:encgradient}
\end{multline}
Next we focus on the boundary term. For $u = \sfU^{-1} v$, we have
\[
	\int_{\p\Omg_\eps} \left(m+\frac{H_1^{\partial\Omega_\eps}}{2}\right)|u|^2 \dd s = \int_{\Sigma_\eps^+} \left(m+\frac{H_1^{\partial\Omega_\eps}}{2}\right)|u|^2 \dd s + \int_{\Sigma_\eps^-} \left(m+\frac{H_1^{\partial\Omega_\eps}}{2}\right)|u|^2 \dd s,
\]
where $\Sg_\eps^\pm$ are given in \eqref{eqn:defsigmapm}. We can express these summands as
\[
	\int_{\Sg_\eps^\pm}\left(m+\frac{H_1^{\partial\Omega_\eps}}{2}\right)|u|^2 \dd s = \int_\Sg 
	\left(m+\frac{H_1^{\partial\Omega_\eps}(\phi_\eps(s,\pm1))}{2}\right)|u(\phi_\eps(s,\pm1))|^2 h^\pm_\eps(s) \dd s
\]
with $h_\eps^\pm(s) = \frac{\varphi_\eps(s,\pm 1)}{\eps}$. Thus, we get
\[
	\int_{\Sg_\eps^\pm}\left(m+\frac{H_1^{\partial\Omega_\eps}}{2}\right)|u|^2 \dd s = \int_{\Sigma}\Big(m + \frac{H_1^{\partial\Omega_\eps}(\phi_\eps(s,\pm1))}2\Big)|v(s,\pm1))|^2 \frac{\varphi_\eps(s,\pm 1)}{\eps} \dd s.
\]
%
Hence, taking into account \eqref{eqn:encgradient}, we can sandwich the quadratic form $b_\eps$ between the forms $b^\pm_\eps$ in the sense that $b_\eps^-[v] \le b_\eps[v]\le b_\eps^+[v]$ for any $v\in\dom b_\eps$ where the forms $b^\pm_\eps$ are defined on the same domain as $b_\eps$ by the expressions
\begin{equation}\label{eq:form_beps}
\begin{aligned}
	b_\eps^\pm[v] &\!=\! (1\!\pm\! c\eps)\int_{\Sg\times(-1,1)}
	|\nb_\Sg v|^2\varphi_\eps \dd s \dd t
	\!+\!
	\frac{1}{\eps^2}\int_{\Sg\tm(-1,1)}
	|\p_t v|^2\varphi_\eps \dd s \dd t\\
	&\quad \!+\!m^2\int_{\Sg\tm(-1,1)}
	|v|^2\varphi_\eps \dd s \dd t
	\!+\!\int_\Sg\left(m\!+\!\frac{H_1^{\partial\Omega_\eps}(\phi_\eps(s,1))}{2}\right)\!
	|v(s,1)|^2\frac{\varphi_
		\eps(s,1)}{\eps}
	\dd s\\
	&\qquad +\int_\Sg\left(m+
	\frac{H_1^{\partial\Omega_\eps}(\phi_\eps(s,-1))}{2}\right)
	|v(s,-1)|^2\frac{\varphi_\eps(s,-1)}{\eps}
	\dd s.
\end{aligned}
\end{equation}

\paragraph{Step 2.}
In order to work in a non-weighted $L^2$-space, we get rid of the term $\varphi_\eps$ in the metric by considering the unitary transform
\[
	\sfV\colon L^2(\Sg\tm(-1,1),\dC^N;\varphi_\eps \dd s \dd t)\arr L^2(\Sg\tm(-1,1),\dC^N;\dd s\dd t),\qquad \sfV v := \sqrt{\varphi_\eps}v,
\]
and define the quadratic forms
\[
	c_\eps^\pm[w] := b_\eps^\pm[\sfV^{-1}w],\qquad \dom c_\eps^\pm :=\sfV(\dom b_\eps).
\]
For $v \in \dom b_\eps^\pm$, we set $w = \sqrt{\varphi_\eps}v$.
By definition of $\varphi_\eps$ we have
\[
	\p_t\varphi_\eps = \sum_{k=1}^{n-1} k\eps^{k+1} t^{k-1} H_k.
\]
By standard rules of differentiation we get
\[	
	\p_t(\varphi_\eps^{-1/2}) 
	= 
	-\frac{	\p_t\varphi_\eps}{2\varphi_\eps^{3/2}}.
\]
Next we compute the expression of $c_\eps^\pm[w]$
\[
\begin{aligned}
	|\p_t(\varphi_\eps^{-1/2} w)|^2 
	&=
	|\p_t(\varphi_\eps^{-1/2})|^2|w|^2 + \frac{1}{\varphi_\eps}|\p_t w|^2 + 
	2(\varphi_\eps^{-1/2})(\p_t(\varphi_\eps^{-1/2}))\Re\Big(\langle\p_t w, w\rangle_{\CC^N}\Big)\\
	&=
	|\p_t(\varphi_\eps^{-1/2})|^2|w|^2 + \frac{1}{\varphi_\eps}|\p_t w|^2 + 
	\varphi_\eps^{-1/2}\p_t(\varphi_\eps^{-1/2})\p_t(|w|^2).
\end{aligned}
\]
Using the auxiliary computation
\[
\begin{aligned}
	\int_{-1}^1\varphi_\eps^{-1/2}\p_t(\varphi_\eps^{-1/2})\p_t(|w|^2)\varphi_\eps \dd t &= -\int_{-1}^1 \frac{\p_t\varphi_\eps}{2\varphi_\eps}\p_t(|w|^2)\dd t\\
	& = \int_{-1}^1\p_t\left(\frac{\p_t\varphi_\eps}{2\varphi_\eps}\right)|w|^2\dd t - \left[\frac{\p_t\varphi_\eps}{2\varphi_\eps}|w|^2\right]_{t=-1}^{t=1}.
\end{aligned}
\]
we find that
\begin{equation}\label{eq:aux_integral}
\begin{aligned}
	\int_{\Sg\tm(-1,1)}|\p_t v|^2\varphi_\eps \dd s \dd t &
	= 
	\int_{\Sg\tm(-1,1)}
	\varphi_\eps|\p_t(\varphi_\eps^{-1/2})|^2|w|^2\dd s\dd t \!+\! \int_{\Sg\tm(-1,1)}|\p_t w|^2 \dd s \dd t\\
	&\quad\quad +\int_{\Sg\tm(-1,1)} \p_t\left(\frac{\p_t\varphi_\eps}{2\varphi_\eps}\right)|w|^2 \dd s \dd t\\
	&\quad\qquad+\int_\Sg \left\{\left[\frac{\p_t\varphi_\eps}{2\varphi_\eps}\right](s,-1)|w(s,-1)|^2
	-\left[\frac{\p_t\varphi_\eps}{2\varphi_\eps}\right](s,1)|w(s,1)|^2
	\right\}\dd s.
\end{aligned}
\end{equation}
Now we pass to the analysis of the next term
\[
\begin{aligned}
	|\nb_\Sg v|^2\varphi_\eps 
	=
	|\nb_\Sg(\varphi_\eps^{-\frac12} w) |^2\varphi_\eps 
	&=
	\varphi_\eps \sum_{j,k=1}^{n-1}g^{j,k} \langle\partial_j(\varphi_\eps^{-\frac12} w),\partial_k(\varphi_\eps^{-\frac12}w)\rangle_{\CC^N} \\
	&= \varphi_\eps \Big(\sum_{j,k = 1}^{n-1} g^{j,k} \partial_j(\varphi_\eps^{-\frac12})\partial_k(\varphi_\eps^{-\frac12})\Big)|w|^2 \\&\qquad + \Big(\sum_{j,k = 1}^{n-1} g^{j,k} \langle\partial_j w,\partial_k w\rangle_{\dC^N}\Big)\\&\qquad + 2\varphi_\eps \sum_{j,k=1}^{n-1}g^{j,k} \Re\Big(\langle \partial_j(\varphi_\eps^{-\frac12})w, \varphi_\eps^{-\frac12} \partial_k w\rangle_{\CC^N}\Big),
\end{aligned}
\]
where we have used the symmetry of $g^{-1}$ to obtain the last terms. Now, remark that for all $x,y \in \CC^{n-1}$ the map $(x,y) \mapsto \langle g^{-1} x,y \rangle_{\CC^{n-1}}$ is a scalar product. Hence for $p \in \{1,\dots,N\}$, if one sets
\[
	x_p = (\partial_j(\varphi_\eps^{-\frac12})w_p)_{j\in\{1,\dots,n-1\}} \in\CC^{n-1},\qquad y_p = (\partial_j w_p)_{j\in\{1,\dots,n-1\}} \in\CC^{n-1},
\]
where we have set $w = (w_1,\dots,w_N)^\top$. There holds
\[
	\sum_{j,k=1}^{n-1}g^{j,k} \langle \partial_j(\varphi_\eps^{-\frac12})w,  \partial_k w\rangle_{\CC^N} = \sum_{p=1}^{N} \langle g^{-1}x_p, y_p\rangle_{\CC^{n-1}}
\]
and using the Cauchy-Schwarz inequality for the scalar product in $\CC^{n-1}$ associated with the metric $g^{-1}$ we obtain
\begin{align*}
	\Big|2\varphi_\eps \sum_{j,k=1}^{n-1}g^{j,k} \Re\Big(\langle \partial_j(\varphi_\eps^{-\frac12})w, \varphi_\eps^{-\frac12} \partial_k w\rangle_{\CC^N}\Big)\Big| &\leq 2\varphi_\eps^{\frac12} \Big|\sum_{j,k=1}^{n-1}g^{j,k} \langle \partial_j(\varphi_\eps^{-\frac12})w, \partial_k w\rangle_{\CC^N}\Big|\\& \leq 2 \varphi_\eps^{\frac12} \sum_{p=1}^N |\langle g^{-1}x_p, y_p\rangle_{\CC^{n-1}}|\\
	& \leq 2 \varphi_\eps^{\frac12} \sum_{p=1}^N (\langle g^{-1}x_p, x_p\rangle_{\CC^{n-1}})^{\frac12}(\langle g^{-1}y_p, y_p\rangle_{\CC^{n-1}})^{\frac12}\\
	& \leq  \varphi_\eps^{\frac12} \sum_{p=1}^N\Big(\eps^{\frac12}\langle g^{-1}y_p, y_p\rangle_{\CC^{n-1}} + \frac1{\eps^{\frac12}}\langle g^{-1}x_p, x_p\rangle_{\CC^{n-1}}\Big).
\end{align*}
But one remarks that we have
\[
\begin{aligned}
	\langle g^{-1}y_p, y_p\rangle_{\CC^{n-1}} &= \sum_{j,k=1}^{n-1}g^{j,k}\partial_j w_p \overline{\partial_k w_p} = |\nabla_\Sigma w_p|^2,\\
	\langle g^{-1}x_p, x_p\rangle_{\CC^{n-1}} &=\sum_{j,k=1}^{n-1}g^{j,k}\partial_j (\varphi_\eps^{-\frac12}) {\partial_k (\varphi_\eps^{-\frac12})}|w_p|^2 = |\nabla_\Sigma (\varphi_\eps^{-\frac12})|^2 |w_p|^2.
\end{aligned}
\]
It gives
\[
	\Big|2\varphi_\eps \sum_{j,k=1}^{n-1}g^{j,k} \Re\Big(\langle \partial_j(\varphi_\eps^{-\frac12})w, \varphi_\eps^{-\frac12} \partial_k w\rangle_{\CC^N}\Big)\Big| \leq \eps^{\frac12}\varphi_\eps^{\frac12}|\nabla_\Sigma w|^2 + \frac{\varphi_\eps^{\frac12}}{\eps^{\frac12}}|\nabla_\Sigma (\varphi_\eps^{-\frac12})|^2 |w|^2.
\]
Using the asymptotics in $L^\infty(\Sigma)$
\[
	\varphi_\eps = \eps +\cO(\eps^2),\qquad	\p_j(\varphi_\eps^{-1/2}) = -\frac{\p_j(\varphi_\eps)}{2\varphi_\eps^{3/2}} = \cO(\eps^{2-3/2}) = \cO(\sqrt{\eps}),\qquad\eps\arr0,
\]
we obtain that there exists some (new) constant $c > 0$ such that for all $\eps$ small enough there holds
\begin{equation}\label{eq:grad_v_est}
\begin{aligned}
	|\nb_\Sg v|^2\varphi_\eps &\le \varphi_\eps|\nb_\Sg(\varphi_\eps^{-1/2})w|^2 + |\nb_\Sg w|^2 + c\eps|\nb_\Sg w|^2 + c\eps|w|^2,\\
		|\nb_\Sg v|^2\varphi_\eps &\geq \varphi_\eps|\nb_\Sg(\varphi_\eps^{-1/2})w|^2 + |\nb_\Sg w|^2 - c\eps|\nb_\Sg w|^2 - c\eps|w|^2.\\
\end{aligned}
\end{equation}
Combining~\eqref{eq:form_beps} with
\eqref{eq:aux_integral},~\eqref{eq:grad_v_est} and the last computation we get that
\[	
c^-_\eps[w] \le c_\eps[v] \le c_\eps^+[w],
\]
where
\begin{equation}\label{key}
\begin{aligned}
	c^\pm_\eps[w] &\!=\! (1\pm c\eps)\int_{\Sg\tm(-1,1)}|\nb_\Sg w|^2\dd s \dd t\\
	&\quad + (1\pm c\eps)\int_{\Sg\tm(-1,1)}\varphi_\eps|\nb_\Sg(\varphi_\eps^{-1/2})|^2|w|^2\dd s \dd t \\
	&\quad+(m^2\pm c\eps)\int_{\Sg\tm(-1,1)}|w|^2\dd s\dd t
	+\frac{1}{\eps^2}\int_{\Sg\tm(-1,1)}|\p_tw|^2\dd s\dd t\\
	&\quad
	+\frac{1}{\eps^2}\int_{\Sg\tm(-1,1)}\varphi_\eps|\p_t(\varphi_\eps^{-1/2})|^2|w|^2\dd s\dd t\\
	&\quad + \frac{1}{\eps^2}\int_{\Sg\tm(-1,1)}\p_t\left(\frac{\p_t\varphi_\eps}{2\varphi_\eps}\right)|w|^2\dd s\dd t\\
	&\quad \!+\! \frac{1}{\eps^2}\int_\Sg
	\left\{\left[\frac{\p_t\varphi_\eps}{2\varphi_\eps}\right](s,\!-1)|w(s,\!-1)|^2
	\!-\!
	\left[\frac{\p_t\varphi_\eps}{2\varphi_\eps}\right](s,\!1)|w(s,\!1)|^2\right\}\dd s\\
	&\quad \!+\!\frac{1}{\eps}
	\int_\Sg\left\{\!\left(m\!+\!\frac{H_1^{\partial\Omega_\eps}(\phi_\eps(s,1))}{2}\right)|w(s,\!1)|^2\! + \!
	\left(m\!+\!\frac{H_1^{\partial\Omega_\eps}(\phi_\eps(s,-1))}{2}\right)|w(s,\!-1)|^2\!\right\}\dd s.\\
\end{aligned}
\end{equation}
\paragraph{Step 3.} Now, we compute and expand in powers of $\eps$ each term of $c_\eps^\pm$ in which appears $\varphi_\eps$ or its derivatives.

We start by dealing with the boundary terms. By Lemma \ref{lem:expcurv}, there holds
\begin{multline*}
	\frac1\eps \int_\Sigma \left(\frac{H_1^{\p\Omega_\eps}(\phi_\eps(s,1))}{2}|w(s,1)|^2 + \frac{H_1^{\p\Omega_\eps}(\phi_\eps(s,-1))}{2}|w(s,-1)|^2\right) \dd s\\= \frac1\eps \int_\Sigma\left( \frac{H_1(s)}2 |w(s,1)|^2 - \frac{H_1(s)}{2}|w(s,-1)|^2\right) \dd s\\\qquad\qquad- \frac12 \int_\Sigma\big(H_1^2(s) - 2 H_2(s)\big)(|w(s,1)|^2 + |w(s,-1)|^2) \dd s \\+ \mathcal{O}(\eps) \int_\Sigma(|w(s,1)|^2 + |w(s,-1)|^2)\dd s.
\end{multline*}
Now, using the representation for $\varphi_\eps$ deduced from Point \eqref{itm:1-metric}~Lemma \eqref{lemma:metric} we find
\[
	\frac{1}{\eps^2}\frac{\p_t\varphi_\eps}{\varphi_\eps} = \frac{1}{\eps^2}\frac{1}{\varphi_\eps}\sum_{k=1}^{n-1} k\eps^{k+1}t^{k-1}H_k = \frac{1}{\eps^2}\frac{\sum_{k=1}^{n-1}k\eps^kt^{k-1}H_k}{1+\sum_{j=1}^{n-1}\eps^jt^jH_j}.
\]
In view of
\[
\begin{aligned}
	\frac{1}{1+\sum_{j=1}^{n-1}\eps^j t^j H_j} &= 1 - \sum_{j=1}^{n-1} \eps^j t^j H_j + \left(\sum_{j=1}^{n-1}\eps^j t^j H_j\right)^2 + \cO(\eps^3)\\
	& = 1-\eps tH_1 -\eps^2t^2H_2 +\eps^2t^2H_1^2+\cO(\eps^3)\\ 
	&=
	1-\eps tH_1 +t^2(H_1^2-H_2)\eps^2 +\cO(\eps^3)
\end{aligned}
\]
we obtain that
\[
\begin{aligned}
	\frac{1}{\eps^2}\frac{\p_t\varphi_\eps}{\varphi_\eps} &= \Big(\sum_{k=1}^{n-1}k\eps^{k-2}t^{k-1}H_k\Big)\left(1-\eps tH_1 +t^2(H_1^2-H_2)\eps^2 +\cO(\eps^3)\right) 
	\\
	&=\frac{H_1}{\eps} - t(H_1^2-2H_2) + \cO(\eps).
\end{aligned}
\]	
Finally, we find that
\[	
	\frac{1}{\eps^2}\frac{\p_t\varphi_\eps(s,\pm1)}{\varphi_\eps} = \frac{H_1(s)}{\eps} \mp (H_1^2(s) -2H_2(s))+\cO(\eps).
\]
It gives
\[
\begin{aligned}
&\frac{1}{\eps^2}\int_\Sigma \left\{\left[\frac{\partial_t \varphi_\eps}{2\varphi_\eps}\right](s,-1)|w(s,-1)|^2 - \left[\frac{\partial_t \varphi_\eps}{2\varphi_\eps}\right](s,1)|w(s,1)|^2\right\} \dd s \\ 
&\qquad\qquad = \frac1\eps \int_\Sigma \left(\frac{H_1(s)}2 |w(s,-1)|^2 - \frac{H_1(s)}2 |w(s,1)|^2\right) \dd s
\\ &\qquad\qquad\qquad\qquad + \frac12 \int_\Sigma (H_1^2(s) - 2H_2(s))(|w(s,-1)|^2 + |w(s,1)|^2)\dd s\\ &\qquad\qquad\qquad\qquad\qquad+\mathcal{O}(\eps) \int_{\Sigma}(|w(s,-1)|^2 + |w(s,1)|^2)\dd s. 
\end{aligned}
\]
We have obtained that
\begin{multline*}
	 \frac{1}{\eps^2}\int_\Sigma \left\{\left[\frac{\partial_t \varphi_\eps}{2\varphi_\eps}\right](s,-1)|w(s,-1)|^2 - \left[\frac{\partial_t \varphi_\eps}{2\varphi_\eps}\right](s,1)|w(s,1)|^2\right\} \dd s  \\ + \frac1\eps \int_\Sigma\left\{ \frac{H_1^{\p\Omega_\eps}(\phi_\eps(s,1))}{2}|w(s,1)|^2 + \frac{H_1^{\p\Omega_\eps}(\phi_\eps(s,-1))}{2}|w(s,-1)|^2\right\} \dd s\\ = \mathcal{O}(\eps) \int_{\Sigma}(|w(s,1)|^2 + |w(s,-1)|^2)\dd s.
\end{multline*}
Next, we expand the term $\frac{1}{\eps^2}\varphi_\eps|\p_t(\varphi_\eps^{-1/2})|^2$. Using the formula
\[
	\p_t(\varphi_\eps^{-1/2}) = -\frac12\frac{\p_t\varphi_\eps}{\varphi_\eps^{3/2}}
\]
we find
\[
	\frac{1}{\eps^2}
	\varphi_\eps|\p_t(\varphi_\eps^{-1/2})|^2=
	\frac{1}{4\eps^2}
	\frac{|\p_t\varphi_\eps|^2}{\varphi_\eps^2}.
\]
Substituting into the above formula the expression
\[
	(\p_t\varphi_\eps)(s,t) = \eps\sum_{k=1}^{n-1} k\eps^kt^{k-1}H_k(s)
\]
we obtain that
\begin{equation}
\begin{aligned}
	\frac{1}{\eps^2}\varphi_\eps|\p_t(\varphi_\eps^{-1/2})|^2 &=
	\frac{1}{4\varphi_\eps^2}\left|\sum_{k=1}^{n-1}k\eps^k t^{k-1} H_k\right|^{2}\\
	&= 
	\frac{1}{4(1+\sum_{k=1}^{n-1}\eps^k t^kH_k)^2}\left|\sum_{k=1}^{n-1}k\eps^{k-1}t^{k-1}H_k\right|^2 = \frac{H_1^2}{4} + \cO(\eps).
\end{aligned}
\label{eqn:pot1}
\end{equation}
Next we consider the term $\frac{1}{\eps^2}\p_t\left(\frac{\p_t\varphi_\eps}{2\varphi_\eps}\right)$.
\[
	\frac{1}{\eps^2}\p_t\left(\frac{\p_t\varphi_\eps}{2\varphi_\eps}\right) = 
	\frac{1}{\eps^2}\left[\frac{\p^2_t\varphi_\eps}{2\varphi_\eps} -\frac{(\p_t\varphi_\eps)^2}{2\varphi_\eps^2}\right].
\]
We expand the two terms on the right-hand side separately. For the first term we obtain
\[
	\frac{1}{\eps^2}\frac{\p_t^2\varphi_\eps}{2\varphi_\eps} =
	\frac{\sum_{k=2}^{n-1}
		k(k-1)\eps^{k-2}t^{k-2}H_k}{2\left(1+\sum_{k=1}^{n-1}\eps^kt^k H_k\right)} = H_2 +\cO(\eps).
\]
For the second term we find
\[
	\frac{1}{\eps^2}\frac{(\p_t\varphi_\eps)^2}{2\varphi_\eps^2} = \frac12\frac{\left(\sum_{k=1}^{n-1}k
		\eps^{k-1} t^{k-1}H_k\right)^2}{\left(1+\sum_{k=1}^{n-1}\eps^k t^k H_k\right)^{2}}  =\frac{H_1^2}{2} + \cO(\eps).
\]
Combining the last two expansions we get
\begin{equation}
	\frac{1}{\eps^2}\p_t\left(\frac{\p_t\varphi_\eps}{2\varphi_\eps}\right) = H_2 - \frac{H_1^2}{2} +\cO(\eps).
	\label{eqn:pot2}
\end{equation}
Finally, we expand the term $\varphi_\eps|\nb_\Sg(\varphi_\eps^{-1/2})|^2$. Observe that there holds
\begin{align*}
	|\nb_\Sg(\varphi_\eps^{-1/2})|^2 &=\sum_{j,k=1}^{n-1}g^{j,k}\p_j(\varphi_\eps^{-1/2})\p_k(\varphi_\eps^{-1/2})\\
	& = \frac1{4\varphi_\eps^3} \sum_{j,k=1}^{n-1}g^{j,k} \partial_j \varphi_\eps \partial_k \varphi_\eps\\ & = \frac1{4\varphi_\eps^3}|\nabla_\Sigma \varphi_\eps|^2.
\end{align*}
Thus, we get
\begin{equation}
	\varphi_\eps|\nabla_\Sigma(\varphi_\eps^{-1/2})|^2= \frac{|\nabla_\Sigma\varphi_\eps|^2}{4\varphi_\eps^2}
	= \frac{\left|\sum_{k=1}^{n-1}\eps^{k} t^k(\nabla_\Sigma H_k)\right|^2}{4\left(1+\sum_{k=1}^{n-1} \eps^k t^k H_k\right)^2} = \cO(\eps^2).
	\label{eqn:pot3}
\end{equation}
Summarizing \eqref{eqn:pot1}, \eqref{eqn:pot2} and \eqref{eqn:pot3}, we have obtained so far that
\[
\begin{aligned}
\varphi_\eps|\nb_\Sg(\varphi_\eps^{-1/2})|^2&= \cO(\eps^2),\\
\frac{1}{\eps^2}\varphi_\eps|\p_t(\varphi_\eps^{-1/2})|^2 &= \frac{H_1^2}{4} + \cO(\eps),\\
\frac{1}{\eps^2}\p_t\left(\frac{\p_t\varphi_\eps}{2\varphi_\eps}\right) &=  H_2- \frac{H_1^2}{2}+\cO(\eps),
\end{aligned}
\]
where the remainders are understood in the $L^\infty(\Sigma\times(-1,1))$-norm.
Using all these expansions we find that there exists a (new) constant $c > 0$ such that for all  $\eps>0$ small enough we can set
\[
\begin{aligned}
c_\eps^\pm[w] &= (1\pm c\eps)\int_{\Sg\tm(-1,1)}|\nb_\Sg w|^2\dd s\dd t + \int_{\Sg\tm(-1,1)} \left(H_2-\frac{H_1^2}{4}\right)|w|^2\dd s \dd t\\
&\quad +\frac{1}{\eps^2}\left[\int_{\Sg\tm(-1,1)}|\p_tw|^2\dd s \dd t+ (m\eps \pm c\eps^3)\int_\Sg(|w(s,-1)|^2+|w(s,1)|^2)\dd s\right]\\
&\qquad\qquad + (m^2\pm c\eps)\int_{\Sg\tm(-1,1)}|w|^2\dd s \dd t
\end{aligned}
\]
such that there holds
\[
	c_\eps^-[w] \leq b_\eps^-[v],\quad b_\eps^{+}[v]\leq c_\eps^+[w].
\]
which gives the claim of Proposition \ref{prop:lbubqf} with $\cU := \sfV\sfU$.
\end{proof}
\section{Proof of the main result}\label{sec:proofmr}
In our setting, a transverse operator appears in the expression of $c_\eps^\pm$ given in Proposition~\ref{prop:lbubqf}. Thanks to Point \eqref{itm:1-1D}~Proposition~\ref{prop:1doperator} we recognize the quadratic form of the square of the transverse Dirac operator $\sfT_{\nu(s)}(\delta)$ introduced in \eqref{eqn:deftransdirac} (for $\delta = m\eps \pm c\eps^{3}$). In particular, its modes depend on the $s$-variable and in \S \ref{subsec:unifest} we give a uniform estimate with respect to the $s$-variable of these modes. \S \ref{subsec:ub} and \ref{subsec:lb} are devoted to the proof of an upper and a lower bound on the eigenvalues of $\scD_\eps^2$, respectively. The proof of Theorem \ref{thm:main} is performed in \S \ref{subsec:proofmain}.
\subsection{A uniform estimate}\label{subsec:unifest}
For further uses, we need the following lemma regarding the modes of the operator
\begin{equation}
	\sfT_s^\delta := \sfT_{\nu(s)}(\delta),
	\label{eqn:optrans-s}
\end{equation}
where $\sfT_{\nu(s)}(\delta)$ is the transverse Dirac operator defined in \eqref{eqn:deftransdirac} and $\delta \geq 0$ is a parameter.  Remark that in this case the modes $\varphi_j^{\delta,\pm} := \varphi_{j,1}^{\delta,\pm}$ defined in Equations \eqref{eqn:posmod1D} and \eqref{eqn:negmod1D} also depend on $s \in \Sigma$.
\begin{lem} There exist constants $C > 0$ and $\delta_0 > 0$ such that for all $j\in\{1,2,\dots,\frac{N}{2}\}$, all $\eta \in \{\pm\}$ and all $\delta \in (0,\delta_0)$ there holds
\[
	\|\varphi_j^{\delta,\eta} - \varphi_j^{\eta}\|_{L^\infty(\Sigma\times(-1,1),\CC^N)}\leq C \delta,\qquad \|\nabla_\Sigma \varphi_j^{\delta,\eta} - \nabla_\Sigma \varphi_j^\eta\|_{L^\infty(\Sigma\times(-1,1),\RR^{n-1}\otimes\CC^N)} \leq  C\delta,
\]
where $\varphi_j^\eta := \varphi_j^{0,\eta}$ corresponds to $\delta = 0$.
\label{lem:asymptexp}
\end{lem}
\begin{proof} Let us fix $j\in \{1,\dots,\frac{N}2\}$. We prove it for $\varphi_j^{\dl,+}$, the proof for $\varphi_j^{\dl,-}$ being similar. Thanks to \eqref{eqn:posmod1D} and \eqref{eqn:modpos1D0}, we remark that for all $t \in [-1,1]$ and $s\in \Sigma$ there holds
\begin{align*}
	\varphi_{j}^{\dl,+}(s,t) - \varphi_j^+(s,t) &= \left(N_\dl k_1(\dl)- \frac12\right) \cos(k_1(\delta)(t+1))\begin{pmatrix}\eps_j\\-i\beta(\nu(s))\eps_j\end{pmatrix} \\&\qquad+ \left(N_\dl E_1(\dl) - \frac12\right)\sin(k_1(\delta)(t+1))\begin{pmatrix}\eps_j\\i\beta(\nu(s))\eps_j\end{pmatrix} \\&\qquad+ N_\dl\delta \sin(k_1(\delta)(t+1))\begin{pmatrix}\eps_j\\-i\beta(\nu(s))\eps_j\end{pmatrix}\\
	&\qquad - \frac12 \left(\cos\left(\frac\pi4(t+1)\right)-\cos(k_1(\delta)(t+1))\right)\begin{pmatrix}\eps_j\\-i\beta(\nu(s))\eps_j\end{pmatrix}\\
	&\qquad - \frac12\left(\sin\left(\frac\pi4(t+1)\right) - \sin(k_1(\delta)(t+1))\right)\begin{pmatrix}\eps_j\\i\beta(\nu(s))\eps_j\end{pmatrix},
\end{align*}
where we have set $N_\dl := N_{\dl,1}$. In particular, for all $t \in [-1,1]$, there holds
\begin{equation}
\begin{aligned}
	&|\varphi_{j}^{\dl,+}(s,t) - \varphi_j^{+}(s,t)| \leq \sqrt{2} \left|N_\dl k_1(\dl)- \frac12\right|  + \sqrt{2} \left|N_\dl E_1(\dl)- \frac12\right| + \sqrt{2} N_\dl\dl\\&\quad\quad \!+\! \frac{\sqrt{2}}2\left|\sin\left(\frac\pi4(t+1)\right) \!-\! \sin(k_1(\delta)(t+1))\right|\! +\! \frac{\sqrt{2}}2\left|\cos\left(\frac\pi4(t+1)\right)\!-\!\cos(k_1(\delta)(t+1))\right|\\&
	\quad \leq \sqrt{2} \left|N_\dl k_1(\dl)- \frac12\right|  + \sqrt{2} \left|N_\dl E_1(\dl)- \frac12\right| + \sqrt{2} N_\dl\dl + \sqrt{2} \left|\frac\pi4-k_1(\delta)\right|, 
	\label{eqn:unifboundphi0}
\end{aligned}
\end{equation}
where we have used the mean-value theorem, that the norm in $\CC^N$ of the vector $\begin{pmatrix}\eps_j\\i\beta(\nu(s))\eps_j\end{pmatrix}$ does not depend on $s$ (because $\beta(\nu)$ is unitary, see \eqref{eqn:betacomm}) and equals $\sqrt{2}$. Similarly, let $(e_1,\dots,e_{n-1})$ be a local coordinate system at $s\in \Sigma$. For $k \in \{1,\dots,{n-1}\}$ there holds
\begin{align*}
	(\partial_k\varphi_{j}^{\dl,+})(s,t) - (\partial_k\varphi_j^{+})(s,t) &= \left(N_\dl k_1(\dl)- \frac12\right) \cos(k_1(\delta)(t+1))\begin{pmatrix}0\\-i\beta(\partial_k\nu(s))\eps_j\end{pmatrix} \\&\qquad+ \left(N_\dl E_1(\dl) - \frac12\right)\sin(k_1(\delta)(t+1))\begin{pmatrix}0\\i\beta(\partial_k \nu(s))\eps_j\end{pmatrix} \\&\qquad+ N_\delta\delta \sin(k_1(\delta)(t+1))\begin{pmatrix}0\\-i\beta(\partial_k \nu(s))\eps_j\end{pmatrix}\\
	&\qquad - \frac12 \left(\cos\left(\frac\pi4(t+1)\right)-\cos(k_1(\delta)(t+1))\right)\begin{pmatrix}0\\-i\beta(\partial_k\nu(s))\eps_j\end{pmatrix}\\
	&\qquad - \frac12
	\left(\sin\left(\frac\pi4(t+1)\right) - \sin(k_1(\delta)(t+1))\right)\begin{pmatrix}0\\i\beta(\partial_k\nu(s))\eps_j\end{pmatrix}.
\end{align*}
It yields
\begin{align*}
	&|\nabla_\Sigma \varphi_{j}^{\delta,+}(s,t) - \nabla_\Sigma \varphi_j^{+}(s,t)|^2 =
	\\
	&\qquad = \sum_{p,k = 1}^{n-1} g^{p,k}
	\big\langle\partial_p \varphi_{j}^{\delta,+}(s,t) - \partial_p \varphi_j^{+}(s,t),\partial_k \varphi_{j}^{\delta,+}(s,t) - \partial_k \varphi_j^{+}(s,t)\big\rangle_{\dC^N}\\
	&\qquad \le \Big(\sum_{p,k=1}^{n-1} g^{p,k}\partial_p \nu(s)\cdot\partial_k\nu(s) \Big)\bigg(\left|N_\dl k_1(\dl)- \frac12\right| + \left|\frac\pi4 - k_1(\delta)\right|+  \left|N_\dl E_1(\dl)- \frac12\right |+ N_\dl\dl\bigg)^2\\
	&\qquad 
	= 
	|\nabla_\Sigma \nu(s)|^2\bigg(\left|N_\dl k_1(\dl)- \frac12\right| + \left|\frac\pi4 - k_1(\delta)\right|+  \left|N_\dl E_1(\dl)- \frac12\right | +N_\dl\dl\bigg)^2\\
	&\qquad = 
	(H_1^2(s) - 2H_2(s))\bigg(\left|N_\dl k_1(\dl)- \frac12\right| + \left|\frac\pi4 - k_1(\delta)\right|+  \left|N_\dl E_1(\dl)- \frac12\right | + N_\dl\dl\bigg)^2,
\end{align*}
where we have used that by Proposition~\ref{prop:alpha_Gamma}
\[
\langle\beta(\p_p\nu(s))\eps_j,
\beta(\p_k\nu(s))\eps_j\rangle_{\dC^{N/2
}}
+
\langle\beta(\p_k\nu(s))\eps_j,
\beta(\p_p\nu(s))\eps_j\rangle_{\dC^{N/2}} = 2\p_p\nu(s)\cdot\p_k\nu(s)
	\rangle
\] and that $|\nabla_\Sigma \nu|^2 = {\rm Tr}(S^2) = H_1^2 - 2H_2$. $\Sigma$ being $C^\infty$, there exists a constant $C >0$ such that $0 \leq H_1^2 - 2H_2 \leq C^2$ and it gives
\begin{equation}
	|\nabla_\Sigma\varphi_{j}^{\dl,+}(s,t) - \nabla_\Sigma\varphi_j^{+}(s,t)| \leq C \left(\left|N_\dl k_1(\dl)- \frac12\right|+  \left|N_\dl E_1(\dl)- \frac12\right| +  N_\dl\dl + \left|\frac\pi4 - k_1(\delta)\right|\right).
	\label{eqn:unifboundderiphi0}
\end{equation}
Hence, the lemma is proved once we prove that the terms on the right-hand side of \eqref{eqn:unifboundphi0} and \eqref{eqn:unifboundderiphi0} are of order $\delta$.

One remarks that there holds
\[
	|\varphi_{j}^{\delta,+}(s,t)|^2 = 2N_\dl^2\Big(\big|k_1(\delta) \cos(k_1(\delta)(t+1)) + \dl \sin(k_1(\delta)(t+1))\big|^2 + E_1(\delta)^2\sin(k_1(\delta)(t+1))^2\Big).
\]
As $\varphi_{j}^{\delta,+}(s,\cdot)$ is normalized in $L^2((-1,1),\CC^N)$, it yields:
\begin{align*}
	1 &= \int_{-1}^1|\varphi_{j}^{\delta,+}(s,t)|^2 \dd t \\&= 2N_\delta^2 \left(\frac14 k_1(\delta)\sin(4k_1(\delta)) + k_1(\delta)^2 + E_1(\delta)^2 +\delta^2 - (E_1(\delta)^2 +\delta^2)\frac{\sin(4k_1(\delta))}{4k_1(\delta)} + \delta \sin(2k_1(\delta))^2\right)\\
	& = 2 N_\dl^2 \left(2 E_1(\delta)^2 -\delta^2 \frac{\sin(4k_1(\delta))}{2k_1(\delta)} +  \delta \sin(2 k_1(\delta))^2\right),
\end{align*}
where we have used that $E_1(\delta)^2 = k_1(\delta)^2 + \delta^2$. Now, using Point \eqref{itm:3-1D}~Proposition \ref{prop:1doperator}, we obtain $1 = N_\dl^2(\frac{\pi^2}{4} + \mathcal{O}(\delta))$ which yields
\[
	N_\delta = \frac{2}{\pi} + \mathcal{O}(\delta).
\]
It gives $E_1(\delta) N_\delta = \frac12 + \mathcal{O}(\delta)$, $k_1(\delta) N_\delta = \frac12 + \mathcal{O}(\delta)$ which, taking into account in \eqref{eqn:unifboundphi0} and \eqref{eqn:unifboundderiphi0}, gives Lemma \ref{lem:asymptexp}.
\end{proof}
\subsection{An upper bound}\label{subsec:ub}
In this paragraph we obtain an upper bound on the eigenvalues of the operator $\scD_\eps^2$ with small $\eps > 0$ in terms of the eigenvalues  of the effective self-adjoint operator $\Upsilon$ associated with the quadratic form in~\eqref{eq:theta}.
The following proposition holds.
\begin{prop}\label{prop:upperbound} 
For all $j \in \NN$, there exists a constant $k> 0$ and $\eps_1 > 0$ such that the inequality
\[
	\mu_j(\scD_\eps^2) \leq \eps^{-2}E_1(m\eps + c\eps^3)^2 + \mu_j(\Upsilon) + k\eps
\]
holds for all $\eps \in (0,\eps_1)$.
\label{prop:ub}
\end{prop}
\begin{proof}
Let us consider a test function of the form
\[
	u(s,t) := \sum_{p = 1}^{\frac{N}2}\big[ f_p^+(s) \varphi_p^{\delta,+}(s,t) + f_p^{-}(s) \varphi_p^{\dl,-}(s,t)\big] 
\]
with $\delta := m\varepsilon + c\varepsilon^3$ where the constant $c$ is as in the definition of $c_\eps^\pm$ in Proposition \ref{prop:lbubqf}. For all $p \in \{1,\dots,\frac{N}2\}$, we assume $f_p^\pm \in H^1(\Sigma)$ and recall that the functions $\varphi_p^{\delta,+}(s,\cdot)$ are the normalized eigenfunctions of $\sfT_s^\delta$ introduced in \eqref{eqn:optrans-s}. We set $f = (f^+,f^-)^\top \in H^1(\Sigma,\CC^N)$ with $(f^\pm)_p = f_p^\pm$ for all $p\in \{1,\dots,\frac{N}2\}$.
One remarks that $u \in \dom c_\eps^+$ and
\begin{multline*}
	c_\eps^+[u] = (1+c\varepsilon) \int_{\Sigma\times(-1,1)} |\nabla_\Sigma u|^{2} \dd s \dd t +\int_{\Sigma} \left(H_2 - \frac{H_1^2}4\right)|f|^2 \dd s \\+ \frac1{\eps^2}\Big(E_1(\delta)^2 - \delta^2\Big)\|f\|_{L^2(\Sigma,\CC^N)}^2 + (m^2 + c\eps)\|f\|_{L^2(\Sigma,\CC^N)}^2,
\end{multline*}
where we have used the
that the functions $\varphi_p^{\delta,\pm}(s,\cdot)$ are orthonormal in $L^2((-1,1),\CC^{N})$  as well as Point (\ref{itm:1-1D}) in~Proposition \ref{prop:1D}. We focus on the term
\begin{multline*}
	\int_{\Sigma\times(-1,1)} |\nabla_\Sigma u|^2 \dd s \dd t = \int_{\Sigma\times (-1,1)} \sum_{j,k = 1}^{n-1} g^{j,k} \langle \partial_j u,\partial_k u\rangle_{\CC^N} \dd s \dd t \\ = \int_{\Sigma\times (-1,1)} \sum_{j,k = 1}^{n-1} \sum_{p,q = 1}^{\frac{N}2} g^{j,k} \langle A_p^{j,\dl} + B_p^{j,\dl}, A_q^{k,\dl} + B_q^{k,\dl}\rangle_{\CC^N} \dd s \dd t\\ = \int_{\Sigma\times (-1,1)} \sum_{j,k = 1}^{n-1} \sum_{p,q = 1}^{\frac{N}2} g^{j,k} \Big(\langle A_p^{j,\delta},A_q^{k,\delta}\rangle_{\CC^N}+ \langle B_p^{j,\delta},B_q^{k,\delta}\rangle_{\CC^N} + \langle A_p^{j,\delta},B_q^{k,\delta}\rangle_{\CC^N} +  \langle B_p^{j,\dl},A_q^{k,\dl}\rangle_{\CC^N}\Big) \dd s \dd t,
\end{multline*}
where for $j \in \{1,\dots, n-1\}$ and $p \in \{1,\dots, \frac{N}2\}$, we have set
\[
	A_p^{j,\dl} = (\partial_j f_p^+) \varphi_p^{\delta,+} + (\partial_j f_p^-) \varphi_p^{\delta,-},\quad B_p^{j,\dl} = f_p^+ \partial_j \varphi_{p}^{\delta,+} + f_p^-\partial_j\varphi_p^{\delta,-}.
\]
One finds
\begin{multline*}
	\langle A_p^{j,\dl},A_q^{k,\dl}\rangle_{\CC^N} = (\partial_j f_p^+) \overline{(\partial_k f_q^+)} \langle\varphi_p^{\delta,+},\varphi_q^{\delta,+}\rangle_{\CC^N} + (\partial_j f_p^+) \overline{(\partial_k f_q^-)} \langle\varphi_p^{\delta,+},\varphi_q^{\delta,-}\rangle_{\CC^N} \\+ (\partial_j f_p^-) \overline{(\partial_k f_q^+)} \langle\varphi_p^{\delta,-},\varphi_q^{\delta,+}\rangle_{\CC^N} + (\partial_j f_p^-) \overline{(\partial_k f_q^-)} \langle\varphi_p^{\delta,-},\varphi_q^{\delta,-}\rangle_{\CC^N}.
\end{multline*}
Using the orthonormality of the $\varphi_j^{\delta,\pm}$ in $L^2((-1,1),\CC^N)$, we obtain
\begin{multline*}
	\int_{\Sigma\times (-1,1)} \sum_{j,k = 1}^{n-1} \sum_{p,q = 1}^{\frac{N}2} g^{j,k} \langle A_p^{j,\dl},A_q^{k,\dl}\rangle_{\CC^N} \dd s \dd t = \int_\Sigma \Big( \sum_{j,k=1}^{n-1}g^{j,k} \sum_{p = 1}^{\frac{N}2} (\partial_j f_p^+) \overline{(\partial_k f_p^+)} + (\partial_j f_p^-) \overline{(\partial_k f_p^-)}\Big) \dd s\\ = \int_{\Sigma} |\nabla_\Sigma f|^2 \dd s.
\end{multline*}
Now, we remark that
\begin{multline*}
	\langle B_p^{j,\delta}, B_q^{k,\delta}\rangle_{\CC^N} = f_p^+ \overline{f_q^+}\langle \partial_j \varphi_p^{\delta,+}, \partial_k \varphi_q^{\delta, +}\rangle_{\CC^N} + f_p^+ \overline{f_q^-}\langle \partial_j \varphi_p^{\delta,+}, \partial_k \varphi_q^{\delta, -}\rangle_{\CC^N} \\+ f_p^- \overline{f_q^+}\langle \partial_j \varphi_p^{\delta,-}, \partial_k \varphi_q^{\delta, +}\rangle_{\CC^N} + f_p^- \overline{f_q^-}\langle \partial_j \varphi_p^{\delta,-}, \partial_k \varphi_q^{\delta, -}\rangle_{\CC^N}.
\end{multline*}
With a slight abuse of notation, we denote
for any $s\in\Sg$ by $\langle\cdot,\cdot\rangle_{\dR^{n-1}\otimes\dC^N}$ the inner product in $T_s\Sg\otimes\dC^N$ with the metric taken into account.
Next, for $p,q \in \{1,\dots,\frac{N}2\}$ and $\eta_1,\eta_2 \in \{\pm\}$, remark that
\begin{align*}
	&\sum_{j,k=1}^{n-1}g^{j,k} \langle\partial_j \varphi_p^{\delta,\eta_1},\partial_k \varphi_q^{\dl,\eta_2}\rangle_{\CC^N}\\
	&\qquad = \langle\nabla_\Sigma \varphi_p^{\delta,\eta_1},\nabla_\Sigma\varphi_q^{\dl,\eta_2}\rangle_{\RR^{n-1}\otimes\CC^N}\\&
	\qquad= \Big\langle\nabla_\Sigma \varphi_p^{\eta_1}\! +\! \Big(\nabla_\Sigma \varphi_p^{\delta,\eta_1}-\nabla_\Sigma \varphi_p^{\eta_1}\Big),\nabla_\Sigma\varphi_q^{\eta_2}\! +\! \Big(\nabla_\Sigma\varphi_q^{\delta,\eta_2}-\nabla_\Sigma\varphi_q^{\eta_2}\Big)\Big\rangle_{\RR^{n-1}\otimes\CC^N}\\
	&\qquad = 
	\langle \nabla_\Sigma \varphi_{p}^{\eta_1},\nabla_\Sigma \varphi_{q}^{\eta_2}\rangle_{\RR^{n-1} \otimes \CC^N} + \langle\nabla_\Sigma \varphi_{p}^{\eta_1},\nabla_\Sigma\varphi_q^{\delta,\eta_2}-\nabla_\Sigma\varphi_q^{\eta_2}\rangle_{\RR^{n-1} \otimes \CC^N} \\&\qquad\qquad\qquad+ \langle\nabla_\Sigma \varphi_p^{\delta,\eta_1}-\nabla_\Sigma\varphi_p^{\eta_1},\nabla_\Sigma \varphi_q^{\eta_2}\rangle_{\RR^{n-1} \otimes \CC^n}\\
	&\qquad\qquad\qquad\quad+
	\langle\nabla_\Sigma \varphi_p^{\delta,\eta_1}-\nabla_\Sigma \varphi_p^{\eta_1},
	\nabla_\Sigma \varphi_q^{\delta,\eta_2}-\nabla_\Sigma \varphi_q^{\eta_2}\rangle_{\dR^{n-1}\otimes\dC^N}.
	\end{align*}
In particular, there holds
\begin{multline*}
	\left|\langle \nabla_\Sigma \varphi_p^{\delta,\eta_1}, \nabla_\Sigma \varphi_q^{\delta,\eta_2}\rangle_{\RR^{n-1}\otimes\CC^N} - \langle \nabla_\Sigma \varphi_p^{\eta_1}, \nabla_\Sigma \varphi_q^{\eta_2}\rangle_{\RR^{n-1}\otimes\CC^N}\right|\\
	\leq 
	\left|\langle\nabla_\Sigma \varphi_{p}^{\eta_1},\nabla_\Sigma\varphi_q^{\delta,\eta_2}-\nabla_\Sigma\varphi_q^{\eta_2}\rangle_{\RR^{n-1} \otimes \CC^N}\right| 
	+ 
	\left|\langle\nabla_\Sigma \varphi_p^{\delta,\eta_1}-\nabla_\Sigma\varphi_p^{\eta_1},\nabla_\Sigma \varphi_q^{\eta_2}\rangle_{\RR^{n-1} \otimes \CC^N}\right|\\ 
	+ \left|\langle \nabla_\Sigma\varphi_p^{\delta,\eta_1}-\nabla_\Sigma\varphi_p^{\eta_1},\nabla_\Sigma\varphi_q^{\delta,\eta_2}-\nabla_\Sigma\varphi_q^{\eta_2}\rangle_{\RR^{n-1}\otimes\CC^N}\right|\\\leq \Big|\nabla_\Sigma \varphi_{p}^{\eta_1}\Big| \Big|\nabla_\Sigma\varphi_q^{\delta,\eta_2}-\nabla_\Sigma\varphi_q^{\eta_2}\Big| + \Big|\nabla_\Sigma \varphi_{p}^{\eta_2}\Big| \Big|\nabla_\Sigma\varphi_q^{\delta,\eta_1}-\nabla_\Sigma\varphi_q^{\eta_1}\Big|\\ + \Big|\nabla_\Sigma\varphi_p^{\delta,\eta_1}-\nabla_\Sigma\varphi_p^{\eta_1}\Big| \Big|\nabla_\Sigma\varphi_q^{\delta,\eta_2}-\nabla_\Sigma\varphi_q^{\eta_2}\Big|\\
	\leq C \Big(\sup_{\eta \in \{\pm\}\atop r \in \{1,\dots,\frac{N}2\}}\|\nabla_\Sigma \varphi_r^{\eta}\|_{L^\infty(\Sigma\times(-1,1), \dR^{n-1}\otimes\CC^N)}\Big)\delta + C^2 \delta^2
	\leq C' \delta,
\end{multline*}
for some new constant $C' > 0$ and $\delta > 0$ small enough. The last inequality, is obtained using Lemma \ref{lem:asymptexp}. For $j \in \{1,\dots,n-1\}$ and $p \in \{1,\dots,\frac{N}2\}$ we set $B_p^j := B_p^{j,0}$ (for $\dl = 0$), we get
\begin{align*}
	&\sum_{j,k = 1}^{n-1}g^{j,k}\Big(\langle B_p^{j,\delta},B_q^{k,\delta}\rangle_{\dC^N}
	 - \langle B_p^{j},B_q^{k}\rangle_{\dC^N}\Big)\\
	&\qquad  = f_p^+\overline{f_q^+}\Big(\langle \nabla_\Sigma \varphi_p^{\delta,+},\nabla_\Sigma \varphi_q^{\dl,+}\rangle_{\RR^{n-1}\otimes\CC^N}-\langle \nabla_\Sigma \varphi_p^{+},\nabla_\Sigma \varphi_q^{+}\rangle_{\RR^{n-1}\otimes\CC^N}\Big) \\&\qquad\qquad+ f_p^+\overline{f_q^-}\Big(\langle \nabla_\Sigma \varphi_p^{\delta,+},\nabla_\Sigma \varphi_q^{\dl,-}\rangle_{\RR^{n-1}\otimes\CC^N} - \langle \nabla_\Sigma \varphi_p^{+},\nabla_\Sigma \varphi_q^{-}\rangle_{\RR^{n-1}\otimes\CC^N}\Big)\\
	& \qquad\qquad +f_p^-\overline{f_q^+}\Big(\langle \nabla_\Sigma \varphi_p^{\delta,-},\nabla_\Sigma \varphi_q^{\dl,+}\rangle_{\RR^{n-1}\otimes\CC^N} - \langle\nabla_\Sigma \varphi_p^{-},\nabla_\Sigma \varphi_q^{+}\rangle_{\RR^{n-1}\otimes\CC^N}\Big) \\&\qquad\qquad+ f_p^-\overline{f_q^-}\Big(\langle \nabla_\Sigma \varphi_p^{\delta,-},\nabla_\Sigma \varphi_q^{\dl,-}\rangle_{\RR^{n-1}\otimes\CC^N} - \langle \nabla_\Sigma \varphi_p^{-},\nabla_\Sigma \varphi_q^{-}\rangle_{\RR^{n-1}\otimes\CC^N}\Big)
\end{align*}
and thus we get
\begin{multline*}
\bigg|\sum_{p,q=1}^{\frac{N}2}\sum_{j,k = 1}^{n-1}g^{j,k}\Big(\langle B_p^{j,\delta},B_q^{k,\delta}\rangle_{\dC^N} - \langle B_p^{j},B_q^{k}\rangle_{\dC^N}\Big)\bigg| \\\leq \sum_{p,q=1}^{\frac{N}2}\sum_{\eta_1,\eta_2\in\{\pm\}}|f_p^{\eta_1}|\cdot|f_q^{\eta_2}|\cdot\left|\langle \nabla_\Sigma \varphi_p^{\delta,\eta_1}, \nabla_\Sigma \varphi_q^{\delta,\eta_2}\rangle_{\RR^{n-1}\otimes\CC^N} - \langle \nabla_\Sigma \varphi_p^{\eta_1}, \nabla_\Sigma \varphi_q^{\eta_2}\rangle_{\RR^{n-1}\otimes\CC^N}\right|\\\leq C'\delta \sum_{p,q=1}^{\frac{N}2}\sum_{\eta_1,\eta_2 \in \{\pm\}}|f_p^{\eta_1}|\cdot|f_q^{\eta_2}|
\leq \frac{C'}2 \delta \sum_{p,q=1}^{\frac{N}2}\sum_{\eta_1,\eta_2\in\{\pm\}}\Big(|f_p^{\eta_1}|^2 + |f_q^{\eta_2}|^2\Big)\leq C'' \delta |f|^2,
\end{multline*}
for some new constant $C'' > 0$ which depends on $N$ but not on $\delta$. Hence, we are left with the understanding of the term $\sum_{p,q=1}^{\frac{N}2}\sum_{j,k = 1}^{n-1}g^{j,k}\langle B_p^{j},B_q^{k}\rangle_{\dC^N}$.
Using \eqref{eqn:modpos1D0} and \eqref{eqn:modneg1D0}, we obtain
\begin{equation}
	\partial_j \varphi_p^+ (s,t) = \frac12 \cos\left(\frac\pi4(t+1)\right)\begin{pmatrix}0 \\ -i\beta(\partial_j \nu)\varepsilon_p\end{pmatrix} + \frac12 \sin\left(\frac\pi4(t+1)\right)\begin{pmatrix}0 \\ i \beta(\partial_j\nu)\varepsilon_p\end{pmatrix}.
	\label{eqn:dpsipos}
\end{equation}
and
\begin{equation}
	\partial_j \varphi_p^- (s,t) = \frac12 \cos\left(\frac\pi4(t+1)\right)\begin{pmatrix}i\beta(\partial_j \nu)^*\varepsilon_p\\0\end{pmatrix} + \frac12 \sin\left(\frac\pi4(t+1)\right)\begin{pmatrix} -i \beta(\partial_j\nu)^*\varepsilon_p\\0\end{pmatrix}.
	\label{eqn:dpsineg}
\end{equation}
Hence, we get
\[
\begin{aligned}
	\int_{-1}^1 \langle\partial_j \varphi_p^+,\partial_k\varphi_q^+\rangle_{\CC^N} \dd t &= \frac14\langle \beta(\partial_j\nu)\varepsilon_p,\beta(\partial_k\nu)\varepsilon_q
	\rangle_{\dC^{\frac{N}{2}}} \int_{-1}^1 \left(\cos\left(\frac\pi4(t+1)\right) - \sin\left(\frac\pi4(t+1)\right)\right)^2\dd t \\
	&= \left(\frac12-\frac1{\pi}\right)\langle \beta(\partial_j\nu)\varepsilon_p,\beta(\partial_k\nu)\varepsilon_q\rangle_{\CC^{\frac{N}2}}
\end{aligned}
\]
as well as
\[
	\int_{-1}^1 \langle\partial_j \varphi_p^+,\partial_k\varphi_q^-\rangle_{\CC^N} \dd t = \int_{-1}^1 \langle\partial_j \varphi_p^-,\partial_k\varphi_q^+\rangle_{\CC^N} \dd t = 0.
\]
Similarly, one obtains
\[
	\int_{-1}^1 \langle\partial_j \varphi_p^-,\partial_k\varphi_q^-\rangle_{\CC^N} \dd t = \left(\frac12 - \frac1\pi\right)\langle\beta(\partial_j \nu)^*\varepsilon_p,\beta(\partial_k\nu)^*\varepsilon_q\rangle_{\CC^{\frac{N}2}}.
\]
Hence, one can write
\begin{multline*}
	\int_{-1}^1\Big(\sum_{j,k = 1}^{n-1}\sum_{p,q = 1}^{\frac{N}2} g^{j,k} \langle B_p^j,B_q^k\rangle_{\CC^N}\Big) \dd t = \sum_{j,k = 1}^{n-1}\sum_{p,q = 1}^{\frac{N}2} g^{j,k} \left(\frac12-\frac1\pi\right)f_p^+ \overline{f_q^{+}}\langle\beta(\partial_j \nu)\varepsilon_p,\beta(\partial_k\nu)\varepsilon_q\rangle_{\CC^{\frac{N}{2}}}\\ + \sum_{j,k=1}^{n-1}\sum_{p,q=1}^{\frac{N}2}g^{j,k}\left(\frac12-\frac1\pi\right)f_p^- \overline{f_q^{-}}\langle\beta(\partial_j \nu)^*\varepsilon_p,\beta(\partial_k\nu)^*\varepsilon_q\rangle_{\CC^{\frac{N}{2}}}\\ = \left(\frac12-\frac1\pi\right)\sum_{p,q = 1}^{\frac{N}2} \big(f_p^+ \overline{f_q^{+}}W_{p,q}+ f_p^- \overline{f_q^{-}}\widetilde{W}_{p,q}\big),
\end{multline*}
where the entries of $\frac{N}{2}\tm\frac{N}{2}$ matrices $W$ and $\wt{W}$ are given by
\[
	W_{p,q} := \left\langle\varepsilon_p , \sum_{j,k = 1}^{n-1} g^{j,k}\beta(\partial_j \nu)^* \beta(\partial_k\nu) \varepsilon_q\right\rangle_{\dC^{\frac{N}{2}}}, \qquad\widetilde{W}_{p,q} = \left\langle\varepsilon_p , \sum_{j,k = 1}^{n-1} g^{j,k}\beta(\partial_j \nu) \beta(\partial_k\nu)^* \varepsilon_q\right\rangle_{\dC^\frac{N}{2}}.
\]
In order to compute $W_{p,q}$ and $\widetilde{W}_{p,q}$ set
\begin{equation}
	V := \sum_{j,k = 1}^{n-1} g^{j,k} \beta(\partial_j \nu)^*\beta(\partial_k \nu),\quad \widetilde{V} := \sum_{j,k = 1}^{n-1} g^{j,k} \beta(\partial_j \nu)\beta(\partial_k \nu)^*
	\label{eqn:defVVtilde}
\end{equation}
One remarks that
\begin{align*}
	\begin{pmatrix} V & 0 \\  0 &\widetilde{V}\end{pmatrix} &= \sum_{j,k=1}^{n-1} g^{j,k}\Gamma(\partial_j\nu) \Gamma(\partial_k\nu) \\&= \frac12 \sum_{j,k=1}^{n-1} g^{j,k}\Big(\Gamma(\partial_j\nu) \Gamma(\partial_k\nu) + \Gamma(\partial_k\nu) \Gamma(\partial_j\nu)\Big) \\&= \sum_{j,k=1}^{n-1} g^{j,k} \partial_j \nu \cdot \partial_k \nu\\& = |\nabla_\Sigma \nu|^2 I_N = (H_1^2 - 2H_2) I_N,
\end{align*}
where we have used that $|\nabla_\Sigma \nu|^2 = {\rm Tr}\,(S^2) =H_1^2 - 2H_2$.
This computation yields that $W_{p,q} = \widetilde{W}_{p,q}= (H_1^2 - 2H_2) \delta_{p,q}$. We obtain
\[
	\int_{\Sigma\times(-1,1)}\Big(\sum_{j,k = 1}^{n-1}\sum_{p,q = 1}^{\frac{N}2} g^{j,k} \langle B_p^j,B_q^k\rangle_{\CC^N}\Big) \dd s\dd t = \left(\frac12 - \frac1\pi\right)\int_{\Sigma}(H_1^2-2H_2)|f|^2 \dd s.
\]
We have just shown that
\[
	\int_{\Sigma\times(-1,1)}\sum_{p,q=1}^{\frac{N}2}\sum_{j,k=1}^{n-1} g^{j,k}\langle B_p^{j,\delta},B_q^{k,\delta}\rangle \dd s \dd t = \left(\frac12 - \frac1\pi\right)\int_{\Sigma}(H_1^2-2H_2)|f|^2 \dd s + \mathcal{O}(\delta) \|f\|_{L^2(\Sigma,\dC^N)}^2.
\]
Let us now focus on the last term and remark that there holds
\[
	\sum_{j,k = 1}^{n-1}\sum_{p,q = 1}^{\frac{N}2} g^{j,k} \langle A_p^{j,\delta}, B_q^{k,\delta}\rangle_{\CC^N} = \sum_{j,k = 1}^{n-1}\sum_{q,p = 1}^{\frac{N}2} g^{k,j} \langle A_q^{k,\delta}, B_p^{j,\delta}\rangle_{\CC^N} = \sum_{j,k = 1}^{n-1}\sum_{q,p = 1}^{\frac{N}2} g^{j,k} \overline{\langle B_p^{j,\delta},A_q^{k,\delta}\rangle_{\CC^N}}
\]
where we have used that the inverse of the metric tensor is symmetric. It yields:
\[
	\sum_{j,k = 1}^{n-1}\sum_{p,q = 1}^{\frac{N}2} g^{j,k} \big(\langle A_p^{j,\delta}, B_q^{k,\delta}\rangle_{\CC^N} + \langle B_p^{j,\delta}, A_q^{k,\delta}\rangle_{\CC^N}\big) = 2 \Re\Big(\sum_{j,k = 1}^{n-1}\sum_{p,q = 1}^{\frac{N}2}g^{j,k}\langle B_p^{j,\delta}, A_q^{k,\delta}\rangle_{\CC^N}\Big).
\]
One has
\begin{multline*}
	\langle B_p^{j,\delta}, A_q^{k,\delta}\rangle_{\CC^N} = f_p^+\overline{(\partial_k f_q^+)} \langle\partial_j\varphi_p^{\delta,+},\varphi_{q}^{\delta,+}\rangle_{\CC^N} + f_p^+ \overline{(\partial_k f_q^-)}\langle \partial_j \varphi_p^{\delta,+},\varphi_q^{\delta,-}\rangle_{\CC^N} \\+ f_p^- \overline{(\partial_k f_q^+)} \langle\partial_j \varphi_p^{\delta,-},\varphi_q^{\delta,+}\rangle_{\CC^N} + f_p^- \overline{(\partial_k f_q^-)}\langle \partial_j \varphi_p^{\delta,-},\varphi_{q}^{\delta,-}\rangle_{\CC^N}.
\end{multline*}
Now for $\eta_1,\eta_2 \in \{\pm\}$ we have
\begin{multline*}
	\Big|\sum_{j,k=1}^{n-1}g^{j,k}\big(\langle \partial_j\varphi_{p}^{\delta,\eta_1},\varphi_q^{\delta,\eta_2}\rangle_{L^2((-1,1),\CC^N)} - \langle \partial_j\varphi_{p}^{\eta_1},\varphi_q^{\eta_2}\rangle_{L^2((-1,1),\CC^N)}\big)\Big| \\= \Big|\sum_{j=1}^{n-1}
	\left(\langle(\nabla_\Sigma \varphi_p^{\delta,\eta_1})_j - (\nabla_\Sigma \varphi_p^{\eta_1})_j, \varphi_q^{\delta,\eta_2}\rangle_{L^2((-1,1),\CC^N)} + \langle(\nabla_\Sigma \varphi_p^{\eta_1})_j, \varphi_q^{\delta,\eta_2}-\varphi_q^{\eta_2}\rangle_{L^2((-1,1),\CC^N)}\right)\Big|\\\leq \sum_{j=1}^{n-1}\|(\nabla_\Sigma\varphi_p^{\delta,\eta_1} - \nabla_\Sigma\varphi_p^{\eta_1})_j\|_{L^2((-1,1),\dC^N)} + \|(\nabla_\Sigma \varphi_p^{\eta_1})_j\|_{L^2((-1,1),\dC^N)}\|\varphi_q^{\delta,\eta_2}-\varphi_q^{\eta_2}\|_{L^2((-1,1),\dC^N)}\\\leq C \delta,
\end{multline*}
where we have used the Cauchy-Schwartz inequality, Lemma \ref{lem:asymptexp}
and that for all $t\in [-1,1]$, $|\nabla_\Sigma \varphi_p^{\eta_1}(s,t)|_{\dR^{n-1}\otimes\CC^N}\leq c $ for some constant $c>0$ uniform in $t\in (-1,1)$, $s \in \Sigma$ and $p\in\{1,\dots,\frac{N}{2}\}$.
Thanks to \eqref{eqn:modpos1D0}, \eqref{eqn:modneg1D0}, \eqref{eqn:dpsipos} and \eqref{eqn:dpsineg},we get
\[
	\int_{-1}^{1}\langle\partial_j\varphi_p^+,\varphi_q^+\rangle_{\dC^N} \dd t  = \left(\frac12-\frac1\pi\right)\langle\beta(\partial_j\nu)\varepsilon_p,\beta(\nu)\varepsilon_q\rangle_{\dC^{\frac{N}{2}}}
\]
as well as
\[
	\int_{-1}^{1}\langle\partial_j\varphi_p^-,\varphi_q^-\rangle_{\dC^N} \dd t  = \left(\frac12-\frac1\pi\right)\langle\beta(\partial_j\nu)^*\varepsilon_p,\beta(\nu)^*\varepsilon_q\rangle_{\dC^{\frac{N}{2}}}.
\]
Moreover, we have
\[
	\int_{-1}^1 \langle\partial_j \varphi_p^+,\varphi_q^-\rangle_{\dC^N} \dd t\! =\! -\frac{i}4
	\langle\beta(\partial_j\nu)\varepsilon_p,\varepsilon_q\rangle_{\dC^{\frac{N}{2}}} \int_{-1}^1 \left(\cos^2\left(\frac{\pi}4(t+1)\right) - \sin^2\left(\frac{\pi}4(t+1)\right)\right)\dd t\! =\! 0
\]
and similarly, one gets
\[
	\int_{-1}^1 \langle\partial_j \varphi_p^-,\varphi_q^+\rangle_{\dC^N} \dd t = 0.
\]
For $j\in \{1,\dots,n-1\}$, if one sets
\[
	\omega_{p,q}^j := i \langle \varepsilon_q,\beta(\partial_j \nu)^*\beta(\nu) \varepsilon_p\rangle_{\dC^{\frac{N}{2}}},\qquad\tilde{\omega}_{p,q}^j := i \langle \varepsilon_q,\beta(\partial_j \nu)\beta(\nu)^* \varepsilon_p\rangle_{\dC^{\frac{N}{2}}}
\]
the matrices $\omega^j := (\omega_{p,q}^j)_{p,q\in\{1,\dots,\frac{N}2\}}$ and $\tilde{\omega}^j := (\tilde{\omega}_{p,q}^j)_{p,q\in\{1,\dots,\frac{N}2\}}$ satisfy
\begin{equation}
(\omega^j)^*= \omega^j,\qquad(\tilde{\omega}^j)^* = \tilde{\omega}^j,
\label{eqn:hermitomega}
\end{equation}
where we have used \eqref{eqn:betacomm} taking into account that $\partial_j \nu \in T\Sigma$. One remarks that
\begin{equation}
\begin{aligned}
	&\int_{-1}^1\Big(\sum_{j,k = 1}^{n-1} \sum_{p,q = 1}^{\frac{N}2}g^{j,k} \langle B_p^{j,\delta}, A_q^{k,\delta}\rangle_{\CC^N}\Big) \dd t\\
	&\qquad = -i \left(\frac12 - \frac1\pi\right)\sum_{j,k = 1}^{n-1} \sum_{p,q = 1}^{\frac{N}2}g^{j,k} \Big(f_p^+\overline{\partial_kf_q^+}\omega_{q,p}^j + f_p^-\overline{\partial_kf_q^-}{\tilde\omega}_{q,p}^j\Big) + |\nabla_\Sigma f|^2 \mathcal{O}(\delta) + |f|^2 \mathcal{O}(\delta)\\
	&\qquad= -i\left(\frac12 - \frac1\pi\right)\sum_{j,k=1}^{n-1} g^{j,k}\left\langle\begin{pmatrix}\omega^j & 0 \\ 0& \tilde{\omega}^j\end{pmatrix} f, \partial_k f\right\rangle_{\CC^N} + |\nabla_\Sigma f|^2 \mathcal{O}(\delta) + |f|^2 \mathcal{O}(\delta).
	\label{eqn:expcrossterm}
\end{aligned}
\end{equation}
In \eqref{eqn:expcrossterm}, the terms of order $\delta$ are obtained using that there exists a constant $C > 0$ such that for all $\eta_1,\eta_2 \in \{\pm\}$ the following uniform estimate holds
\[
	\left|\sum_{p,q = 1}^{\frac{N}2} \sum_{j,k=1}^{n-1} g^{j,k} f_p^{\eta_1}\overline{(\partial_k f_q^{\eta_2})}\right| \leq C \left(|\nabla_\Sigma f^{\eta_2}|^2 +|f^{\eta_1}|^2\right).
\]
We are lead to define the matrix-valued $1$-form
\[
	\Omega := \sum_{j=1}^{n-1}\Omega^j\dd s_j,\qquad \Omega^j := \begin{pmatrix}\omega^j & 0\\ 0 & \tilde{\omega}^j\end{pmatrix} \in \CC^{N\times N},
\]
which acts in local coordinates as $(\Omega f)_k = \sum_{j=1}^{n-1} g^{j,k}\Omega^j f$ for $k\in\{1,2,\dots,n-1\}$ for $f\in H^1(\Sigma,\dC^N)$. 
Now remark that \eqref{eqn:expcrossterm} becomes:
\begin{multline*}
	\!\!\!\!\!\int_{-1}^1\Big(\sum_{j,k = 1}^{n-1} \sum_{p,q = 1}^{\frac{N}2}g^{j,k} \langle B_p^{j,\delta}, A_q^{k,\delta}\rangle_{\dC^N}\Big) \dd t \!=\!- i \left(\frac12-\frac1\pi\right)\sum_{j,k = 1}^{n-1}g^{j,k}\big\langle \Omega^j f, \partial_k f\big\rangle_{\dC^N}
	+ |\nabla_\Sigma f|^2 \mathcal{O}(\delta) + |f|^2 \mathcal{O}(\delta).
\end{multline*}
Now, consider the quantity
\[
\begin{aligned}
	&\int_{\Sigma} \Big|\nabla_\Sigma f -i \left(\frac12-\frac1\pi\right)\Omega f\Big|^2 \dd s = \int_\Sigma |\nabla_\Sigma f|^2 \dd s - 2\Re\left(i\int_\Sigma \sum_{j,k=1}^{n-1} g^{j,k} \left\langle \left(\frac12-\frac1\pi\right)\Omega^j f, \partial_k f\right\rangle_{\dC^N} \dd s\right)\\ 
	&\qquad\qquad\qquad\qquad\qquad\qquad\qquad\qquad\qquad+ \Big(\frac12-\frac1\pi\Big)^2 \int_\Sigma \Big(\sum_{j,k=1}^{n-1}g^{j,k}\langle \Omega^j f,\Omega^k f\rangle_{\dC^N}\Big)\dd s.
\end{aligned}	
\]
Remark that
\[
	\sum_{j,k = 1}^{n-1}g^{j,k} \langle\Omega^j f, \Omega^k f\rangle_{\dC^N}
	 = \sum_{j,k = 1}^{n-1}g^{j,k} \langle \Omega^k\Omega^j f, f\rangle_{\dC^N},
\]
where we have used that $(\Omega^k)^* = \Omega^k$ thanks to \eqref{eqn:hermitomega}.
To compute this term, let us introduce the potential
\[
	W := \sum_{j,k = 1}^{n-1}g^{j,k} \Omega^k \Omega^j = \begin{pmatrix} w & 0 \\ 0& \widetilde{w}\end{pmatrix}.
\]
where we have set $w := \sum_{j,k = 1}^{n-1}g^{j,k} \omega^k \omega^j$ and $\widetilde{w} := \sum_{j,k = 1}^{n-1}g^{j,k} \tilde{\omega}^k \tilde{\omega}^j$.
Remark that for all $p,q\in \{1,\dots, \frac{N}2\}$, there holds
\begin{align*}
w_{p,q} = \sum_{j,k=1}^{n-1}g^{j,k} \sum_{r=1}^{\frac{N}2}\omega_{p,r}^k\omega_{r,q}^j &= -\sum_{j,k=1}^{n-1}g^{j,k} \sum_{r=1}^{\frac{N}2}\langle \varepsilon_r,\beta(\partial_k\nu)^*\beta(\nu) \varepsilon_p\rangle_{\CC^{\frac{N}2}}\langle\varepsilon_q,\beta(\partial_j\nu)^*\beta(\nu)\varepsilon_r\rangle_{\CC^{\frac{N}2}}\\
& = - \sum_{j,k = 1}^{n-1} g^{j,k} \sum_{r= 1}^{\frac{N}2} \overline{\langle\beta(\partial_k\nu)^*\beta(\nu) \varepsilon_p,\varepsilon_r\rangle_{\CC^{\frac{N}2}}} {\langle\beta(\nu)^*\beta(\partial_j\nu)\varepsilon_q,\varepsilon_r\rangle_{\CC^{\frac{N}2}}}\\
& = -\sum_{j,k=1}^{n-1}g^{j,k} \langle\beta(\nu)^*\beta(\partial_j\nu)\varepsilon_q,\beta(\partial_k\nu)^*\beta(\nu)\varepsilon_p\rangle_{\CC^{\frac{N}2}}\\
& = -\sum_{j,k=1}^{n-1}g^{j,k} \langle\varepsilon_q,\beta(\partial_j\nu)^*\beta(\nu)\beta(\partial_k\nu)^*\beta(\nu)\varepsilon_p\rangle_{\CC^{\frac{N}2}}\\
& =  \sum_{j,k=1}^{n-1}g^{j,k} \langle\varepsilon_q,\beta(\partial_j\nu)^*\beta(\partial_k\nu)\varepsilon_p\rangle_{\CC^{\frac{N}2}}\\
& = \langle \varepsilon_q, V \varepsilon_p\rangle_{\CC^{\frac{N}2}}\\
& = (H_1^2 - 2H_2) \delta_{p,q}
\end{align*}
where ${V}$ is the scalar potential defined in \eqref{eqn:defVVtilde} and where we have used that $\beta(\nu)\beta(\nu)^* = I_{\frac{N}2}$, that $\partial_k \nu \cdot \nu =0$ and \eqref{eqn:betacomm}. Similarly, one proves that $\widetilde{w}_{p,q} = (H_1^2 - 2H_2) \delta_{p,q}$.

In consequence, we obtain for a constant $c'>0$:
\begin{multline}
	c_\eps^+[u] \leq (1+c'\varepsilon) \int_\Sigma\Big(\Big|\nabla_\Sigma f -i \left(\frac12-\frac1\pi\right) \Omega f\Big|^2 +\Big( \Big(\frac12 + \frac2{\pi^2}\Big)H_2-\frac{H_1^2}{\pi^2}\big) |f|^2\Big)\dd s \\+ \frac1{\varepsilon^2}\big(E_1(m\varepsilon+c\varepsilon^3)+ c'\varepsilon^3 \big) \int_\Sigma |f|^2 \dd s.
	\label{eqn:upboundfin}
\end{multline}
It remains to prove that as defined, for all $j \in \{1,\dots,n-1\}$ there holds $\Omega^j =  -i \Gamma(\partial_j\nu)\Gamma(\nu)$. Remark that
\[
	- i \Gamma(\partial_j\nu)\Gamma(\nu) = \begin{pmatrix}-i\beta(\partial_j\nu)^*\beta(\nu)&0\\0&-i\beta(\partial_j\nu)\beta(\nu)^*\end{pmatrix}.
\]
Remark that the matrix $M^j = (m_{p,q}^j)_{p,q\in\{1,\dots,\frac{N}2\}}$ of the linear map $ \CC^{\frac{N}2}\ni X\mapsto-i\beta(\partial_j\nu)^*\beta(\nu)X$ in the orthonormal basis $(\varepsilon_1,\dots,\varepsilon_{\frac{N}2})$ verifies for all $p,q \in \{1,\dots,\frac{N}2\}$:
\begin{align*}
	m_{p,q}^j 	&= \langle -i\beta(\partial_j\nu)^*\beta(\nu)\varepsilon_q,\varepsilon_p\rangle_{\CC^{\frac{N}2}}\\
			& = -i\langle\varepsilon_q,\beta(\nu)^*\beta(\partial_j\nu) \varepsilon_p\rangle_{\CC^{\frac{N}2}}\\& = i \langle\varepsilon_q,\beta(\partial_j\nu)^*\beta(\nu)\varepsilon_p\rangle_{\CC^{\frac{N}2}} = \omega^j_{p,q}.
		\end{align*}
Similarly, one would prove the matrix $\tilde{M}^j = (\tilde{m}_{p,q}^j)_{p,q\in\{1,\dots,\frac{N}2\}}$ of the linear map $ \CC^{\frac{N}2}\ni X\mapsto-i\beta(\partial_j\nu)\beta(\nu)^*X$ in the orthonormal basis $(\varepsilon_1,\dots,\varepsilon_{\frac{N}2})$ verifies $\tilde{M}^j = \tilde{\omega}^j$ which gives
\[
	\Omega^j = -i \Gamma(\partial_j\nu)\Gamma(\nu).
\]
Hence, \eqref{eqn:upboundfin} writes
\[
	c_\eps^+[u] \leq (1+ c'\eps)\mathfrak{u}[f] + \frac1{\eps^2}\big(E_1(m\eps + c\eps^3) + c' \eps^3\big) \|f\|_{L^2(\Sigma,\CC^N)}^2
\]
and as $\|u\|_{L^2(\Sigma\times(-1,1),\CC^N)} = \|f\|_{L^2(\Sigma,\CC^N)}$, the min-max principle gives Proposition \ref{prop:upperbound}.
\end{proof}
\subsection{A lower bound}\label{subsec:lb}
In this subsection we obtain a lower bound on the eigenvalues of the operator $\scD_\eps^2$ with small $\eps > 0$ in terms of the eigenvalues  of the effective self-adjoint operator $\Upsilon$ associated with the quadratic form in~\eqref{eq:theta}.
Now let us recall the construction from~\cite[Section 4.5]{HOP18}, for which it is useful to employ the identification
\[
	L^2(\Sg\times(-1,1),\dC^N) \simeq L^2(\Sg,\cG),\qquad
	\cG := L^2((-1,1),\dC^N).
\]
Recall that for any Banach space $B$ the gradient $\nabla_\Sg\colon C^1(\Sigma,B)\arr C^0(T\Sigma, B)$ acts in local coordinates of $\Sg$ as
\begin{equation}\label{eq:gradient_local}
	(\nabla_\Sg u)_j = \sum_{k=1}^{n-1} g^{j,k}\p_k u,\qquad j\in\{1,2,\dots,n-1\}. 
\end{equation}
Let us fix an arbitrary point $s_0\in\Sg$ and introduce $C^\infty$-smooth maps $\Theta,\Theta^*\colon \Sg\arr\cB(\cG)$ by
\[
	\Theta(s) = \Theta_{ \nu(s), \nu(s_0) }
	\qquad\text{and}\qquad
	\Theta^*(s) = \Theta_{ \nu(s), \nu(s_0) }^*,
\]
where the unitary map $\Theta_{x,y}\colon\cG\arr\cG$ is as in Point \eqref{itm:5-unitary}~Proposition~\ref{prop:1doperator}
 and $\nu(s)$, $\nu(s_0)$ are the unit normal vectors to the manifold $\Sg$ at the points $s$ and $s_0$, respectively.
Hence, there exists a constant $C > 0$ such that for every $u\in C^0(\Sg,\cG)$ at every point $s\in\Sg$ there holds
\begin{equation}\label{eq:Theta_bounds}
	\|(\nb_\Sg\Theta)(s)u(s)\|_{T_s\Sg\otimes\cG}\le C\|u(s)\|_\cG,\qquad
	\|(\nb_\Sg\Theta^*)(s)u(s)\|_{T_s\Sg\otimes\cG}\le C\|u(s)\|_\cG.
\end{equation}
Throughout this subsection we use the shorthand notation $\sfT_s^\dl := \sfT_{\nu(s)}(\dl)$ 
as in \eqref{eqn:optrans-s} with $\dl = m\eps-c\eps^3$, where the constant $c$ is as in definition of the quadratic form $c_\eps^-$ given in Proposition~\ref{prop:quadratic_forms}.
Furthermore, let $\pi(s)$ be the orthogonal projector in the Hilbert space $L^2((-1,1),\dC^N)$ onto the $N$-dimensional subspace  $\ker((\sfT_s^\dl)^2 - E_1^2(\dl))$ (see Points \eqref{itm:2-1D} and \eqref{itm:4-1D}~Proposition~\ref{prop:1doperator}). 
Denote by $\Pi$ the orthogonal projector in $L^2(\Sg\times(-1,1),\dC^N)$ given by 
\[
	\big(\Pi u\big)(s,t) := \pi(s)u(s,\cdot)(t),
\]
and set $\Pi^\bot := I - \Pi$. Both $\Pi$ and $\Pi^\bot$ define in the canonical way bounded operators in $L^2(T\Sg)\otimes L^2((-1,1), \dC^N)$, to be denoted by the same symbols.
With all the above notation we formulate and prove an auxiliary lemma
on the commutator between $\nb_\Sg$ and the projector $\Pi$, which is reminiscent of~\cite[Lemma 4.16]{HOP18}.
\begin{lem}\label{lem:commutator}
	For $s \in \{0,1\}$, the map $[\nb_\Sg,\Pi]u := \nb_\Sg(\Pi u) - \Pi(\nb_\Sg u)$ defined for $u\in C^1(\Sg)\otimes L^2((-1,1),\dC^N)$ extends by density to a bounded operator
	\[
		[\nb_\Sg,\Pi]\colon L^2(\Sg)\otimes L^2((-1,1),\dC^N)\arr L^2(T\Sg)\otimes L^2((-1,1),\dC^N).
	\]
Moreover, if $u \in H^1(\Sg)\otimes L^2((-1,1),\dC^N)$ then $[\nb_\Sg,\Pi] u \in H^1(T\Sg)\otimes L^2((-1,1),\dC^N)$. The same conclusion holds for $[\nb_\Sg,\Pi^\bot] \equiv -[\nb_\Sg,\Pi]$.
\end{lem}
\begin{proof}
	Let us set $\pi_0 := \pi(s_0)$.
	By Point \eqref{itm:5-unitary}~Proposition~\ref{prop:1doperator}, one has the representation
	$\pi(s) = \Theta(s)\pi_0\Theta(s)^*$. As $\pi_0$ does not depend on $s$, we can perform
	a direct computation in local coordinates to get
	that for any $s\in\Sg$ and any $u\in C^1(\Sg)\otimes L^2((-1,1),\dC^N)$ one has for all $j \in \{1,\dots,n-1\}$
	\begin{align}\label{eq:commutator}
		([\nb_\Sg,\Pi]u)_j(s)
		&= \Big(\sum_{k=1}^{n-1}g^{j,k} (\partial_k \Theta)(s)\pi_0\Theta(s)^*u(s) + \Theta(s) \pi_0 (\partial_k \Theta^*)(s) u(s)\Big)\\& = (\nabla_\Sigma \Theta)_j(s)\pi_0\Theta(s)^*u(s) + \Theta(s) \pi_0 (\nabla_\Sigma \Theta^*)_j(s) u(s).\nonumber
	\end{align}
	Using~\eqref{eq:Theta_bounds} we estimate for any fixed $s\in\Sg$
	\[
	\begin{aligned}
	\|(\nb_\Sg\Th)(s)\pi_0\Th^*(s)u(s)\|_{T_s\Sg\otimes\cG}
	&\le C\|\pi_0\Th^*(s)u(s)\|_\cG\\ &\le C\|\pi_0\|_{\cB(\cG)}\|\Th^*(s)\|_{\cB(\cG)}\|u(s)\|_\cG \le C\|u(s)\|_\cG,\\[0.6ex]
	\|\Th(s)\pi_0(\nb_\Sg\Th^*)(s)u(s)\|_{T_s\Sg\otimes\cG}
	&\le \|\Th(s)\|_{\cB(T_s\Sg\otimes\cG)}\|\pi_0\|_{\cB(\cG)}\|(\nb_\Sg\Th^*)(s)u(s)\|_{T_s\Sg\otimes\cG}\le C\|u(s)\|_\cG,\\
	\end{aligned}
	\]
	where we used that $\Th(s), \Th^*(s)$ are unitary
	and that $\pi_0$ is an orthogonal projector. Moreover, we get
	using~\eqref{eq:commutator} and the triangle inequality for the norm that
	\[
	\begin{aligned}
		\|[\nb_\Sg,\Pi]u\|^2_{L^2(T\Sg)\otimes\cG} &=
		\int_\Sg \big\|(\nb_\Sg\Theta)(s)\pi_0\Theta^*(s) u(s) + \Theta(s)\pi_0
		(\nabla_\Sg\Theta^*)(s)u(s)\|^2_{T_s\Sg\otimes\cG}\dd s\\
		& \le
		2\int_\Sg \big\|(\nb_\Sg\Theta)(s)\pi_0\Theta^*(s) u(s)\big\|^2_{T_s\Sg\otimes\cG}\dd s\\
		&\qquad\qquad +
		2\int_\Sg\big\| \Theta(s)\pi_0
		(\nabla_\Sg\Theta^*)(s)u(s)\|^2_{T_s\Sg\otimes\cG}\dd s \\
		&\le 4C^2\int_\Sg\|u(s)\|_\cG^2\dd s = 
		4C^2\|u\|^2_{L^2(\Sg)\otimes\cG}.
	\end{aligned}
	\]
	Hence, the continuity result follows. To prove the mapping properties of $[\nb_\Sg ,\Pi]$ between the Sobolev spaces of order 1, it
	is enough to remark that~\eqref{eq:commutator} is continuously differentiable with respect to $s$ because $\Th$ and $\Th^*$ are $C^\infty$-smooth.
\end{proof}
Now we have all the tools to obtain a lower bound on the eigenvalues of $\scD_\eps^2$.
\begin{prop}
For all $j \in \NN$, there exists a constant $k> 0$ and $\eps_1 > 0$ such that the inequality
	\[
	\mu_j(\scD_\eps^2) \geq \eps^{-2}E_1(m\eps - c\eps^3)^2 + \mu_j(\Upsilon) - k\eps
	\]
	\label{prop:lb}
holds for all $\eps \in (0,\eps_1)$.	
\end{prop}
\begin{proof}
	Recall that the closed, densely defined, symmetric, and semi-bounded quadratic form $c_\eps^-$ in the Hilbert space $L^2(\Sg\times(-1,1), 
	\dC^N,\dd s\dd t)$  is defined in Proposition~\ref{prop:quadratic_forms} and according to that proposition
	combined with the min-max principle
	one has $\mu_j(c_\eps^-) \le \mu_j(\scD_\eps^2)$ for all $j\in\dN$.
	Set $\dl := m\eps-c\eps^3$, where the constant $c$ is as in definition of the quadratic form $c_\eps^-$. Let $v\in\dom c_\eps^-$ be arbitrary. Using the spectral theorem and Proposition~\ref{prop:1doperator} we obtain the following lower bound
	\begin{equation}\label{eq:1destimate}
	\begin{aligned}
		&\int_{\Sg\times(-1,1)}|\p_t v|^2\dd s\dd t
		+ \dl\int_\Sg\left(|v(s,-1)|^2 + |v(s,1)|^2\right)\dd s\\
		&\qquad \ge (E_1^2(\dl) -\dl^2)\|\Pi v\|^2_{L^2(\Sg\times(-1,1),\dC^N)} + (E_2^2(\dl)-
		\dl^2)\|\Pi^\bot v\|^2_{L^2(\Sg\times(-1,1),\dC^N)},
	\end{aligned}
	\end{equation}
	where the eigenvalues $E_p(\dl)$, $p\in\dN$, are defined as in Point \eqref{itm:2-1D}~Proposition \ref{prop:1doperator}.
	It follows from Point \eqref{itm:2-1D}~Proposition~\ref{prop:1doperator}
	that
	\begin{equation}\label{eq:1stEv}
	E_1^2(\dl) \leq \dl^2 + \frac{\pi^2}4
	\end{equation}
	and that
	\begin{equation}\label{eq:2ndEv}
		E_2^2(\dl) \ge \dl^2 + \frac{9\pi^2}{16}.
	\end{equation}
	Using the pointwise orthogonality $\langle\Pi v(s,\cdot),\Pi^{\bot}v(s,\cdot)\rangle_{L^2((-1,1),\dC^N)} = 0$, $s\in\Sg$, one gets
	\[
	\begin{aligned}
		&\int_{\Sg\times(-1,1)}\left(H_2 -\frac{H_1^2}{4} + m^2-c\eps\right)|v|^2\dd s \dd t\\
		&\qquad =
		\int_{\Sg\times(-1,1)}\left(H_2 -\frac{H_1^2}{4} + m^2-c\eps\right)|\Pi v|^2\dd s \dd t\\ &\qquad\qquad\qquad+\int_{\Sg\times(-1,1)}\left(H_2 -\frac{H_1^2}{4} + m^2-c\eps\right)|\Pi^\bot v|^2\dd s \dd t.
	\end{aligned}
	\]
	The estimate~\eqref{eq:1destimate} combined with~\eqref{eq:2ndEv} and with the above identity yields
	\begin{equation}\label{eq:cepslowerbound}
	\begin{aligned}
		c_\eps^-[v] &\ge (1-c\eps)\int_{\Sg\times(-1,1)}
		|\nb_\Sg v|^2\dd s\dd t\\
		&\qquad + \left\langle\Pi v,\left(H_2-\frac{H_1^2}{4}\right)\Pi v\right\rangle_{L^2(\Sg\times(-1,1),\dC^N)}\\
		&\qquad + 
		\frac{E_1(\dl)^2-\dl^2}{\eps^2} \|\Pi v\|^2_{L^2(\Sg\times(-1,1),\dC^N)} +(m^2- c\eps)\|\Pi v\|_{L^2(\Sigma\times(-1,1),\CC^N)}^2
		\\
		&\qquad
		+ \left\langle\Pi^{\bot} v, \left(H_2-\frac{H_1^2}{4}\right)\Pi^\bot v\right\rangle_{L^2(\Sg\times(-1,1),\dC^N)}\\
		&\qquad+ 
		\left(\frac{9\pi^2}{16\eps^2} +m^2 - c\eps\right)
		\|\Pi^\bot v\|^2_{L^2(\Sg\times(-1,1),\dC^N)}.
	\end{aligned}
	\end{equation}
	Now, we would like to analyse in more details the first term on the right-hand side. Notice that the smoothness of the mappings $\Th$ and $\Th^*$ yields that $\Pi v,\Pi^\bot v\in H^1(\Sg, L^2((-1,1),\dC^N))$.
	 In the following the norm $\|\cdot\|$ and the inner product $\langle\cdot,\cdot\rangle$ are taken either in the Hilbert space
	 $L^2(\Sg, L^2((-1,1),\dC^N))$ or in the Hilbert space  $L^2(T\Sg, L^2((-1,1),\dC^N))$ as no confusion can arise.
	 Then we get
	\begin{equation}\label{eq:gradient_expansion}
		\|\nb_\Sg v\|^2 = \|\nb_\Sg(\Pi v)\|^2 + \|\nb_\Sg(\Pi^\bot v)\|^2 + 2\Re\langle\nb_\Sg(\Pi v),\nb_\Sg(\Pi^\bot v)\rangle,
	\end{equation}
	and
	\[
	\begin{aligned}
		\langle\nb_\Sg(\Pi v), \nb_\Sg(\Pi^\bot v)\rangle 
		& =
		\langle\nb_\Sg(\Pi\Pi v), \nb_\Sg(\Pi^\bot\Pi^\bot v)\rangle\\
		&=
		\langle([\nb_\Sg,\Pi] + \Pi\nb_\Sg)\Pi v, 
		([\nb_\Sg,\Pi^\bot] + \Pi^\bot\nb_\Sg)\Pi^\bot v\rangle\\
		& = 
		\langle[\nb_\Sg,\Pi]\Pi v, 
		[\nb_\Sg,\Pi^\bot]\Pi^\bot v\rangle
		+
		\langle\Pi\nb_\Sg\Pi v, 
		[\nb_\Sg,\Pi^\bot]\Pi^\bot v\rangle\\
		&\qquad+
		\langle[\nb_\Sg,\Pi]\Pi v, 
		\Pi^\bot\nb_\Sg\Pi^\bot v\rangle
		+
		\langle \Pi\nb_\Sg\Pi v, 
		\Pi^\bot\nb_\Sg\Pi^\bot v\rangle\\
		&=: J_1+J_2+J_3+J_4.
	\end{aligned}
	\]
	Relying on the definition of $\Pi$ and $\Pi^\bot$ we immediately conclude that $J_4 = 0$. By Lemma~\ref{lem:commutator} we estimate with $c_0,c_1 >0$ independent of $\eps$
	\[
	\begin{aligned}
		|J_1| &\le c_0\|\Pi v\|\cdot\|\Pi^\bot v\| \le \frac12c_0\eps\|\Pi v\|^2 + \frac{c_0}{2\eps}\|\Pi^\bot v\|^2,\\
		|J_2| &\le c_1 \|\nb_\Sg\Pi v\|\cdot\|\Pi^\bot v\| \le \frac12 c_1\eps\|\nb_\Sg\Pi v\|^2 + \frac{c_1}{2\eps}\|\Pi^\bot v\|^2.
	\end{aligned}
	\]
	Next, using self-adjointness of the projector $\Pi^\bot$ and that by Lemma~\ref{lem:commutator} one has $\Pi^\bot[\nb_\Sg,\Pi]\Pi v\in H^1(T\Sg,L^2((-1,1),\dC^N))$ and thanks to \cite[Prop. 2.2]{Taylor}
	(see also~\cite[\S I.2]{Chavel}) we can perform an integration by parts to get
	\[
	\begin{aligned}
		J_3 &= \langle\Pi^\bot[\nb_\Sg,\Pi]\Pi v,\nb_\Sg\Pi^\bot v\rangle\\
		&= \int_\Sg
		(\Pi^\bot[\nb_\Sg,\Pi]\Pi v,\nb_\Sg\Pi^\bot v)_{T_s\Sigma\otimes L^2((-1,1),\dC^N)}\dd s\\
		&=
		-\int_\Sg
		((\div_\Sg(\Pi^\bot[\nb_\Sg,\Pi]\Pi v))(s),(\Pi^\bot v)(s))_{ L^2((-1,1),\dC^N)}\dd s\\
		&=
		-\langle\div_\Sg(\Pi^\bot[\nb_\Sg,\Pi]\Pi v),\Pi^\bot v\rangle,
 	\end{aligned}
	\]
	which gives the following upper bound
	\[
		|J_3| \le \|\div_\Sg(\Pi^\bot[\nb_\Sg,\Pi]\Pi v)\|\ \|\Pi^\bot v\|.
	\]
	Recall that by \cite[\S I.1, Eq. (26)]{Chavel} in the local coordinates on $\Sg$ for an arbitrary vector field $A = (A_j)$ in the Sobolev space $H^1(T\Sg)$ one has
	\begin{equation}\label{eq:divergence}
		\div_\Sg A = \sum_{j=1}^{n-1}\left(\p_j A_j +\sum_{k=1}^{n-1}\G_{kj}^j A_k\right),
	\end{equation}
	where $\G_{ij}^k$, $i,j,k\in\{1,2,\dots, n-1\}$, are the  \emph{Christoffel symbols}, which are determined by the identity $\nb_{e_j}e_i = \sum_{k=1}^{n-1} \G_{ij}^k e_k$. In our setting the $j$-th component of the vector $\Pi^\bot[\nb_\Sg,\Pi]\Pi v$ is given in  \eqref{eq:commutator} by
	\[
		(\Pi^\bot[\nb_\Sg,\Pi]\Pi v)_j = 
		\Th\pi_0^\bot\Th^* \sum_{k=1}^{n-1} g^{j,k}\left((\p_k\Th)\pi_0\Th^* + \Th\pi_0\p_k\Th^*\right)(\Pi v),
	\]
where we omit the dependence on $s\in\Sg$ for the sake of brevity. Furthermore, the projector $\pi_0$ does not depend on $s$ while $\Th$ is a $C^\infty$-smooth map and does not depend on $\eps$. Therefore, using formula~\eqref{eq:divergence} and boundedness of Christoffel symbols on
a smooth manifold $\Sg$ we get that there exists a constant $c_2 > 0$
such that the  following estimate holds
\[
\begin{aligned}
	&\|\div_\Sg(\Pi^\bot[\nb_\Sg,\Pi]\Pi v)\|_{L^2(\Sg\times(-1,1),\dC^N)}\\
	&\qquad\le \sqrt{c_2}\left(\|\Pi v\|_{L^2(\Sg\times(-1,1),\dC^N)} + 
	\|\nb_\Sg(\Pi v)\|_{L^2(T\Sg,L^2((-1,1),\dC^N))}\right),
\end{aligned} 
\]
which gives
\[
\begin{aligned}
	|J_3| &\le  c_2\eps\big(\|\Pi v\|^2_{L^2(\Sg\times(-1,1),\dC^N)} + \|\nb_\Sg(\Pi v)\|^2_{L^2(T\Sg, L^2((-1,1),\dC^N))}\big) + \frac{1}{2\eps}
	\|\Pi^\bot v\|^2_{L^2(\Sg\times(-1,1),\dC^N)}.
\end{aligned}\]
Hence, from~\eqref{eq:gradient_expansion}
with a suitable $c_3 > 0$
\[
\begin{aligned}
	\|\nb_\Sg v\|^2_{L^2(T\Sg, L^2((-1,1),\dC^N))} & \ge \|\nb_\Sg(\Pi v)\|^2_{L^2(T\Sg, L^2((-1,1),\dC^N))}\\
	&\qquad -2|\langle\nb_\Sg(\Pi v), \nb_\Sg(\Pi^\bot v)\rangle_{L^2(T\Sg, L^2((-1,1),\dC^N))}|\\
	&\ge (1-c_3\eps)\|\nb_\Sg(\Pi v)\|^2_{L^2(T\Sg, L^2((-1,1),\dC^N))} - c_3\eps\|\Pi v\|^2_{L^2(\Sg,L^2((-1,1),\dC^N))}\\ &\qquad-\frac{c_3}{\eps}\|\Pi^\bot v\|^2_{L^2(\Sg,L^2((-1,1),\dC^N))}.
\end{aligned}
\] 
Next, we substitute all the obtained estimates into~\eqref{eq:cepslowerbound}. One notices that a lower bound on all the terms involving $\Pi^\bot v$ is given
by the expression
\[
\begin{aligned}
	&\left(\frac{9\pi^2}{16\eps^2}+m^2-c\eps -\frac{c_3(1-c\eps)}{\eps}\right)\|\Pi^\bot v\|^2_{L^2(\Sg\times(-1,1),\dC^N)}\\
	&\qquad\qquad + \left\langle\Pi^\bot v,\left(H_2-\frac{H_1^2}{4}\right)\Pi^\bot v\right\rangle_{L^2(\Sg\times(-1,1),\dC^N)},
\end{aligned}
\]
which can be bounded from below by $\frac{\pi^2}{2\eps^2}\|\Pi^\bot v\|^2_{L^2(\Sg\times(-1,1),\dC^N)}$ for all $\eps \in (0,\eps_0)$ with sufficiently small $\eps_0 > 0$.
For $\eps \in (0,\eps_0)$ we can simply estimate
\begin{equation}\label{eq:lowerbound}
	c_\eps^-[v] \ge d_{\eps}[\Pi v] +\frac{\pi^2}{2\eps^2}\|\Pi^\bot v\|^2_{L^2(\Sg\times(-1,1),\dC^N)}.
\end{equation}
where $d_{\eps}$ is the quadratic form in the Hilbert space $\ran\Pi$ defined on $\Pi(\dom c_\eps^-)$ by the expression
\[
\begin{aligned}
	d_{\eps}[u] & = (1-c_3\eps)\int_{\Sg\times(-1,1)}
	|\nb_\Sg u|^2\dd s \dd t+\left\langle u,\left(H_2-\frac{H_1^2}{4}\right) u\right\rangle_{L^2(\Sg\times(-1,1),\dC^N)} \\
	&\qquad + 
	\left(\frac{E_1(\dl)^2-\dl^2}{\eps^2} +m^2 - c_4\eps\right) \|u\|^2_{L^2(\Sg\times(-1,1),\dC^N)},
\end{aligned}\]
with a suitable constant $c_4 > 0$.
It follows from inequality~\eqref{eq:lowerbound} and the min-max principle  that for $
\eps\in (0,\eps_0)$
\begin{equation}\label{eq:bndmin}
	\mu_j(c_\eps^-) \ge j-\text{th smallest element of the set } \left\{\mu_j(d_\eps)\right\}_{j\in\dN} \bigcup \left\{\frac{\pi^2}{2\eps^2}\right\}
\end{equation}
It remains to show a connection between $d_\eps$ and the quadratic form $\mathfrak{u}$ for the effective operator introduced in~\eqref{eq:theta}. Clearly, a generic element in $u\in\ran\Pi$ can be  represented in a unique way as
\begin{equation}\label{eq:u}	
	u = \sum_{p=1}^{N/2}\big[f_p^+(s)\varphi_p^{\dl,+}(s,t) + f_p^-(s)\varphi^{\dl,-}_p(s,t)\big],
\end{equation}
where $\dl = m\eps - c\eps^3$, for all $p\in \{1,\dots,N/2\}$, $f_p^\pm\in L^2(\Sg)$, and $\varphi_p^{\dl,\pm}$ are  orthonormal eigenfunctions of $\sfT_s^\dl$ corresponding to the eigenvalues $\pm E_1(\dl)$  (introduced in~\eqref{eqn:optrans-s}).
Note also that $u\in\Pi(\dom c_\eps^-)$ 
is equivalent to the fact that
in the decomposition~\eqref{eq:u} one has $f_p^\pm\in H^1(\Sg)$ for all $p\in\{1,2,\dots,N/2\}$. 
Using the unitary map defined via decomposition~\eqref{eq:u} as
\[
V\colon \ran\Pi\arr L^2(\Sg,\dC^N),\qquad Vu := f = \{f_1^+,\dots f_{N/2}^+, f_1^-,\dots,f_{N/2}^-\}
\]
 we get repeating the same computation as in the last step of the proof of Proposition~\ref{prop:upperbound} that
the quadratic form $d_\eps$ is unitarily equivalent by $V$ to a form $\wt{d}_\eps$ defined on $H^1(\Sg,\dC^N)$ which verifies for some constant $c_5 > 0$:
\[
	\begin{aligned}
	\wt{d}_\eps[u] \geq (1-c_{5}\eps)\mathfrak{u}[f] + \left(\frac{E_1(\dl)^2-\dl^2}{\eps^2}+ m^2- c_{5}\eps\right)\|f\|^2_{L^2(\Sg,\dC^N)},
	\end{aligned}
\]
defined in the Hilbert space $L^2(\Sg,\dC^N)$.
For any fixed $j\in\dN$, 
we get using the estimate~\eqref{eq:1stEv}, the inequality~\eqref{eq:bndmin},
the expression for the quadratic form $\wt{d}_\eps$, and the min-max principle that there exists $\eps_1 > 0$ such that for any $\eps \in (0,\eps_1)$ one has
\[
	\begin{aligned}
	\mu_j(\scD_\eps^2)&\ge\mu_j(c_\eps^-) \ge \min\left\{\mu_j(d_\eps),\frac{\pi^2}{2\eps^2}\right\}\\
	& = \mu_j(d_\eps) = \mu_j(\wt{d}_\eps) \geq
	\frac{E_1(m\eps-c\eps^3)^2}{\eps^2}+ \left(m^2 - \frac{\dl^2}{\eps^2}\right) + \mu_j(\Upsilon) - c_6\eps,
	\end{aligned}
\]
with a constant $c_6 >0$.
As $m^2 - \frac{\dl^2}{\eps^2} = \cO(\eps)$ we obtain the expected result.
\subsection{Proof of Theorem~\ref{thm:main}}\label{subsec:proofmain}
To prove Theorem \ref{thm:main} it is enough to use Point \eqref{itm:3-1D}~Proposition~ \ref{prop:1doperator} to note that for all $\delta > 0$ we have
\[
	E_1(\delta)^2 = \delta^2 + k_1(\delta)^2 = \frac{\pi^2}{16} + \delta - \frac{4}{\pi^2} \delta^2 + \delta^2 +\cO(\delta^3),\qquad\delta\arr0.
\]
In particular, there holds
\[
	E_1(m\eps \pm c\eps^3)^2 = \frac{\pi^2}{16} + m\eps - \frac{4}{\pi^2} m^2\eps^2 + m^2 \eps^2 + \cO(\eps^3),\qquad\eps\arr0,
\]
and combining it with Propositions \ref{prop:upperbound} and \ref{prop:lb} we get that for all $j \in \N$ there holds
\[
	\mu_j(\scD_\eps^2) = \frac{\pi^2}{16\eps^2} + \frac{m}{\eps} - \frac{4}{\pi^2} m^2 + m^2 +  \mu_j(\Upsilon) + \cO(\eps),
	\qquad\eps\arr0,
\]
which is precisely Theorem \ref{thm:main}.

\subsection*{Acknowledgements}VL gratefully acknowledges the support by the grant No. 21-07129S of the Czech Science Foundation.

\begin{appendix}
	\section{The effective operator in two dimensions}\label{app}
	The aim of this appendix is to derive a more explicit formula for the quadratic form $\mathfrak{u}$ in~\eqref{eq:theta} of the effective operator $\Upsilon$ in the case of two dimensions ($n =2 $). This more explicit representation will help us to compare our results with~\cite{K23}.
	
	First of all, we notice that $N = 2^{\lfloor\frac{3}{2}\rfloor} = 2$.
	In this setting, $(\aa_j)_{j=1}^3$ are given by the standard Pauli matrices
	\[
		\aa_1 = \begin{pmatrix} 0 & 1\\
		1&0\end{pmatrix},\qquad\aa_2=\begin{pmatrix} 0 &-i\\
		i &0\end{pmatrix},\qquad \aa_3 = \begin{pmatrix} 1 &0\\
		0&-1\end{pmatrix}.
	\]
	Let $L$ be the length of $\Sg$
	and assume that $\Sg$ is parametrized by the unit-speed mapping $\s\colon [0,L]\arr\dR^2$, $|\dot\s|\equiv 1$ in the clockwise direction. Let $\mathbb{T} = \RR/\!\raisebox{-.65ex}{\ensuremath{L\ZZ}}$ and let $\kp\colon[0,L]\arr\dR$ stand for the curvature of $\Sg$ with the convention that $\kp$ is non-positive if the curve $\Sg$ surrounds a convex domain.
	We introduce the notation $\nu(s) = \nu(\s(s))$ for the normal vector to $\Sg$ at the point $\s(s)$ pointing outwards of the domain surrounded by $\Sg$.
	Using the Frenet formula we get $\nu'_1(s) = -\kp(s)\nu_2(s) $ and $\nu_2'(s) = \kp(s)\nu_1(s)$. 
	The one-form $\Omg = \Omg^1\dd s$ defined in~\eqref{eq:Omega} becomes a $2\times2$ matrix-valued function on $\Sg$, which is expressed in the
	arc-length parametrization of $\Sg$ by
	\[
	\begin{aligned}
		\Omg^1(s) 
		& = -i\Gamma(\nu'(s))\Gamma(\nu(s)) = 
		-i\big(\aa_1\nu_1'(s)+\aa_2\nu_2'(s)\big)\big(\aa_1\nu_1(s)+\aa_2\nu_2(s)\big)
		\\
		& =  -i\begin{pmatrix}
		0 & -\kp(s)(\nu_2(s)+i\nu_1(s))\\
		-\kp(s)(\nu_2(s) -i\nu_1(s))&0\end{pmatrix}
		\begin{pmatrix}
		0 & \nu_1(s) -i\nu_2(s)\\
		\nu_1(s)+i\nu_2(s) & 0\end{pmatrix}\\
	&=-\kp(s)\aa_3. 		
	\end{aligned} 
	\]
	Hence, we derive from~\eqref{eq:theta} using that
	in the two-dimensional setting $H_1^2 = \kp^2$ and $H_2 = 0$
	\[
	\begin{aligned}
	\mathfrak{u}[f] & =\int_0^L\left|
	f'(s) + i\left(\frac{1}{2}-\frac{1}{\pi}\right)\aa_3\kp(s)f(s)\right|^2\dd s -\int_0^L\frac{\kp^2(s)}{\pi^2}|f(s)|^2\dd s,\qquad \dom\mathfrak{u}= H^1(\Sg,\dC^2),
	\end{aligned}
	\]
	which coincides with the quadratic form for the effective operator in~\cite{K23}.
Actually, this effective operator can be simplified. Indeed, if one considers the following unitary map $W \colon L^2(\mathbb{T}) \to L^2(\mathbb{T})$ such that for $f= (f^+,f^-)^\top \in L^2(\mathbb{T})$ there holds
\[
	Wf = \begin{pmatrix}e^{-iV}f_+\\e^{iV}f_-\end{pmatrix},\quad \text{where for } s\in \mathbb{T} \text{ we have set } V(s) = \left(\frac{1}{2}-\frac{1}{\pi}\right)\Big(\int_0^s\kappa(\eta)\dd\eta + \frac{2\pi}Ls\Big)
\]
and for $f =(f^+,f^-)^\top\in W^*( H^1(\mathbb{T}))$ there holds
\[
\mathfrak{u}[Wf] = \mathfrak{u}_{\rm mag}[f_+] + \mathfrak{u}_{\rm mag}[f_-],
\]
where
\[
	\mathfrak{u}_{\rm mag}[g] = \int_{0}^L\Big(\left|g'(s) - i\frac{(\pi-2)}Lg(s) \right|^2  - \frac{\kappa^2}{\pi^2}|g(s)|^2\Big)\dd s,\quad \dom{\mathfrak{u}_{\rm mag}} = W^*( H^1(\mathbb{T})) = H^1(\mathbb{T}).
\]
Remark that the periodic boundary conditions are preserved through the unitary map $W$. Indeed, by definition $V(0) = 0$ and the total curvature identity as well as the clockwise parametrization of $\Sigma$ give $V(L) = 0$ because
\[
	\int_0^L\kappa(s)\dd s = - 2\pi.
\]
$\mathfrak{u}_{\rm mag}$ is associated with the unique self-adjoint operator $\Upsilon_{\rm mag}$ acting in $L^2(\mathbb{T})$ like
\[
	\Upsilon_{\rm mag} = \Big(-i\frac{\der}{\der s} + \frac{(\pi-2)}L\Big)^2 - \frac{\kappa^2}{\pi^2}
\]
on the domain $H^2(\dT)$. 
We have shown that the effective operator $\Upsilon$ is unitarily
equivalent to the orthogonal sum $\Upsilon_{\rm mag}\oplus\Upsilon_{\rm mag}$. 
The operator $\Upsilon_{\rm mag}$ can be understood as a one-dimensional magnetic Schr\"odinger operator with electric potential induced by the curvature. The magnetism here is produced by the non-simply connectedness of the initial domain $\Omega_\varepsilon$.
\end{appendix}

\end{proof}

\end{document}

%% file: commands.tex
\newcommand\nb{\nabla}
\newcommand\tm{\times}

\newcommand{\beq}{\begin{equation} \begin{split}}
\newcommand{\eeq}{\end{split} \end{equation}}
\newcommand\Sg{\Sigma}
\newcommand\Omg{\Omega}

\newcommand\scD{\mathscr{D}}
\makeatletter
\def\section{\@startsection{section}{1}\z@{.9\linespacing\@plus\linespacing}%
	{.7\linespacing} {\fontsize{13}{14}\selectfont\bfseries\centering}}
\def\paragraph{\@startsection{paragraph}{4}%
	\z@{0.3em}{-.5em}%
	{$\bullet$ \ \normalfont\itshape}}

\@addtoreset{equation}{section}
\makeatother

\renewcommand\and{\qquad\text{and}\qquad}

\newcommand\dl{\delta}

\newcommand{\comm}[1]{}

\def\bm1{\mathbbm{1}}
\def\G{\Gamma}

\def\s{\sigma}

\def\p{\partial}

\def\Re{{\rm Re}\,}

\def\arr{\rightarrow}

\def\aa{\alpha}

\def\s{\sigma}

\def\p{\partial}

\def\kp{\kappa}

\newcommand\dd{{\,\mathrm{d}}}
\newcommand\der{{\mathrm{d}}}

\def\sfU{\mathsf{U}}
\def\sfV{\mathsf{V}}

\newcounter{counter_a}

\usepackage[latin1]{inputenc}
\usepackage[T1]{fontenc}

\numberwithin{figure}{section}
\numberwithin{equation}{section}
\theoremstyle{plain}
\newtheorem*{thm*}{Theorem}
\newtheorem{thm}{Theorem}[section]

\newtheorem{lem}[thm]{Lemma}
\newtheorem{prop}[thm]{Proposition}

\theoremstyle{remark}

\theoremstyle{plain}

\newcommand{\beu}{\begin{equation*}}
\newcommand{\eeu}{\end{equation*}}
\newcommand{\besu}{\begin{equation*}
\begin{aligned}}
\newcommand{\eesu}{\end{aligned}
\end{equation*}}
\newcommand{\bes}{\begin{equation}
\begin{aligned}}
\newcommand{\ees}{\end{aligned}
\end{equation}}

\newcommand\cB{\mathcal B}

\newcommand\cG{\mathcal G}
\newcommand\cH{\mathcal H}

\newcommand\CC{\mathbb C}
\newcommand\NN{\mathbb N}
\newcommand\RR{\mathbb R}
\newcommand\ZZ{\mathbb Z}

\newcommand\eps{\varepsilon}

\newcommand\wt{\widetilde}

\newcommand\void[1]{}

\def\eps{\varepsilon}
\def\ran{{\rm ran\,}}

      \def\dC{{\mathbb C}}

   \def\dN{{\mathbb N}}   
      \def\dR{{\mathbb R}}
\def\dS{{\mathbb S}}   \def\dT{{\mathbb T}}

   \def\sfT{{\mathsf T}}   \def\sfU{{\mathsf U}}
\def\sfV{{\mathsf V}}

   \def\cB{{\mathcal B}}   
      
\def\cG{{\mathcal G}}   \def\cH{{\mathcal H}}   
      
      \def\cO{{\mathcal O}}
      
      \def\cU{{\mathcal U}}
      \def\cX{{\mathcal X}}
\def\cY{{\mathcal Y}}

\def\N{\mathbb{N}}

\renewcommand{\div}{\mathrm{div}\,}

\newcommand{\dom}{\mathrm{dom}\,}